\newif\ifTR		

\TRtrue

\ifTR
\documentclass[a4paper,11pt]{article}
\usepackage{authblk}
\usepackage{natbib}
\else
\documentclass[authoryear,review]{elsarticle}

\fi
\usepackage[cm]{fullpage}
\usepackage{caption}
\usepackage{subcaption}
\usepackage{amsthm}

\newtheorem{proposition}{Proposition}
\newtheorem{theorem}{Theorem}

\newtheorem{corollary}{Corollary}

\usepackage{amsmath,amsfonts,amssymb}
\usepackage{graphicx}
\usepackage{paralist}

\usepackage{comment}
\usepackage{bm}
\usepackage{url}
\usepackage{graphics}
\usepackage[capitalize,noabbrev]{cleveref}
\usepackage{multirow,longtable}
\usepackage{rotating}
\usepackage[dvipsnames]{xcolor}
\usepackage{adjustbox}
\usepackage{xspace}
\usepackage{booktabs}
\usepackage{pdflscape}
\usepackage{tikz}
\usepackage{float}
\usepackage{amsmath,amsfonts,amssymb, bm}
\usepackage[normalem]{ulem}
\usepackage{nicefrac}
\usepackage{color} 
\usepackage{pgfplots}
\usepackage{booktabs}
\usepackage{pifont}
\newcommand*\rot{\rotatebox{90}}
\newcommand*\OK{\ding{51}}
\pgfplotsset{compat = newest}
\usepackage{booktabs} 
\usepackage{subcaption}
\usepackage{algorithm}
\usepackage{algpseudocode}
\usepackage[T1]{fontenc}
\usepackage{lmodern}
\usepackage{tabularx}
\usepackage{graphicx}

\usetikzlibrary{shapes,arrows}

\allowdisplaybreaks

\newcommand{\argmin}{\mathrm{argmin}}

\newcommand\probname{Crowdshipping Assignment Problem with Compensation-Driven Acceptance Behavior\xspace}
\newcommand\probabbr{CAPCAB\xspace}
\newcommand\probnonsepabbr{CAPCABwS\xspace}

\newcommand\pnlty{\rho\xspace}
\newcommand\utlty{\mu\xspace}
\newcommand\odno{|J|\xspace}

\newcommand{\rev}[1]{{\color{black}#1}}

\newcommand{\R}{\mathbb{R}}
\ifTR
\fi
\newcommand\crossmark[1][]{%
	\tikz[scale=0.07,#1]{
		\fill(0,0)--(0.1,0) .. controls (0.5,0.4) .. (1,0.7)--(0.9,0.7) ..  controls (0.5,0.5) ..(0,0.1) --cycle;
		\fill(1,0.1)--(0.9,0.1) .. controls (0.5,0.3) .. (0,0.7)--(0.1,0.7) .. controls (0.5,0.4) ..(1,0.2) --cycle;
	}%
}

\newcommand{\appendixproof}[2]{
		\begin{proof}[\emph{\bfseries{Proof of #1.}}] 
		#2
		\end{proof}
}

\tikzstyle{decision} = [diamond, draw, text width=4.5em, text badly centered, node distance=3cm, inner sep=0pt]
\tikzstyle{block} = [rectangle, draw, text width=5em, text centered, rounded corners, minimum height=4em]
\tikzstyle{line} = [draw, -latex']

\newcommand{\myabstract} 
{
	High demand, rising customer expectations, and government regulations are forcing companies to increase the efficiency and sustainability of urban (last-mile) distribution. Consequently,  several new delivery concepts have been proposed that increase flexibility for customers and other stakeholders. One of these innovations is crowdsourced delivery, where deliveries are made by occasional drivers who wish to utilize their surplus resources (unused transport capacity) by making deliveries in exchange for some compensation. The potential benefits of crowdsourced delivery include reduced delivery costs and increased flexibility (by scaling delivery capacity up and down as needed). The use of occasional drivers poses new challenges because (unlike traditional couriers) neither their availability nor their behavior in accepting delivery offers is certain. The relationship between the compensation offered to occasional drivers and the probability that they will accept a task has been largely neglected in the scientific literature. Therefore, we consider a setting in which compensation-dependent acceptance probabilities are explicitly considered in the process of assigning delivery tasks to occasional drivers. We propose a mixed-integer nonlinear model that minimizes the expected delivery costs while identifying optimal assignments of tasks to a mix of professional and occasional drivers and their compensation. We propose an exact two-stage solution algorithm that allows to decompose compensation and assignment decisions for generic acceptance probability functions and show that the runtime of this algorithm is polynomial under mild conditions. Finally, we also study a more general case of the considered problem setting, show that it is NP-hard and propose an approximate linearization scheme of our mixed-integer nonlinear model. 		
	The results of our computational study show clear advantages of our new approach over existing ones. They also indicate that these advantages remain in dynamic settings when tasks and drivers are revealed over time and in which case our method constitutes a fast, yet powerful heuristic.
}

\title{The role of individual compensation and acceptance decisions in crowdsourced delivery}	

\ifTR

\author{Alim Bu\u{g}ra \c{C}{\i}nar}
\author{Wout Dullaert}
\author{Markus Leitner}
\author{Rosario Paradiso}
\author{Stefan Waldherr}
\affil{Department of Operations Analytics, Vrije Universiteit Amsterdam, The Netherlands. 
	\texttt{\{a.b.cinar|m.leitner|wout.dullaert|r.paradiso|s.m.g.waldherr\}@vu.nl}}

\else

\fi

\begin{document}	

\ifTR
\maketitle

\begin{abstract}
	\myabstract 
	
\end{abstract}

\fi

\section{Introduction} \label{sec:intro}

Fueled by the Covid-19 pandemic, e-tailing and same-day home delivery have continued to grow. This increase in demand and rising customer expectations are forcing providers to improve their efficiency in last-mile distribution. At the same time, technological advances and the proliferation of the sharing economy, such as Airbnb, car-sharing systems, and ride-hailing systems like Uber, have popularized the gig economy. In this new paradigm, instead of hiring employees on long-term contracts, employers (companies or individuals) can offer small, short-term tasks (usually through an online platform) and compensation for performing the task. Gig workers can then agree to perform that task in exchange for the advertised compensation. A prime example of the gig economy is crowdsourced delivery \citep{savelsbergh2022challenges}, which is one of the most prominent research directions in city logistics \citep{kaspi2022directions}. Instead of using a traditional logistics system, a company can choose to offer delivery tasks to regular customers or other independent couriers (often referred to as occasional drivers) and compensate them accordingly. Crowdsourced delivery which has the potential to reduce overall delivery costs and increase the flexibility of delivery capacity, has already attracted interest from both industry and academia. \citet{alnaggar2021crowdsourced} and \citet{savelsbergh2022challenges} provide comprehensive overviews and classifications of industry applications and scientific contributions. 

Crowdsourced delivery adds significant complexity to the planning process because the behavior and availability of occasional drivers is often not known in advance. Integrating the behavior of occasional drivers into planning decisions, especially with respect to task acceptance, is a major challenge in this domain \citep{savelsbergh2022challenges}. In the case of instant delivery (e.g., food orders), the inability to find an occasional driver for an offered task may lead to customer dissatisfaction or significant additional costs. In the classical models, occasional drivers are assumed to accept all tasks for which the compensation is high enough to compensate for the additional costs incurred by the drivers for the detour, see e.g. 
\citet{Archetti2016, arslan_crowdsourced_2019}. More recently, acceptance probabilities have been investigated through behavioral studies and surveys of occasional drivers, see e.g. \citet{devari_crowdsourcing_2017,Le2019}.
Such acceptance probabilities were used as fixed parameters in the task allocation problem by \citet{gdowska2018stochastic} and \citet{santini2022probabilistic}.
The effect of compensation schemes on acceptance behavior was investigated by \citet{dayarian_crowdshipping_2020}, but only through a sensitivity analysis of fixed compensation options and acceptance thresholds drawn from a normal distribution. 
\rev{\citet{le2021designing} propose a compensation scheme that defines occasional drivers' behavior as threshold constraints. Thresholds are identified through survey data analysis and any offer above these thresholds is assumed to be accepted. }
\citet{Barbosa2022} and \citet{hou_optimization_2022} integrate the impact of compensation on acceptance probabilities within task allocation. 
However, \citet{Barbosa2022} assume that the detour of all occasional drivers is the same for each task. Consequently, the compensations of all occasional drivers are identical. \citet{hou_optimization_2022} propose a sequential approach in which the compensation values are optimized for fixed task allocations that are decided in the first step. 
\rev{
Our study differs from \citet{le2021designing} by using probability functions to describe occasional driver behavior instead of threshold constraints; from \citet{Barbosa2022} by considering the characteristics of each occasional driver-task pair rather than assuming identical values for all; and from \citet{hou_optimization_2022} by evaluating all possible occasional driver task combinations and compensation values to find a global optimal solution, rather than using a sequential approach without a guarantee of optimality.}

\paragraph{Contribution and outline}
This paper is the first to incorporate task-acceptance probabilities of occasional drivers, which depend on the operator's compensation decisions, into an exact solution framework. To this end, we introduce a mixed-integer nonlinear programming (MINLP) formulation that optimizes the decision process by allowing for full flexibility in the compensations offered and the associated acceptance probabilities of occasional drivers. To model the latter, we introduce acceptance probability functions that estimate the probability that an occasional driver will accept a task. This estimation can be based on historical information, the attributes of occasional drivers and tasks, external factors, and the compensations offered.
Our model is applicable to several established crowdsourced delivery settings, see, e.g., \citet{sampaio_chapter_2019} for a classification. These include in-store deliveries in which occasional drivers (e.g., regular customers) are initially located at the pick-up location and multi-depot deliveries in which they first need to travel to a pick-up location as occurring, e.g., in meal delivery.
In our model, referred to as the \probname (\probabbr), we aim to minimize the total expected cost of assigning tasks by outsourcing them to either occasional (crowdsourced) drivers or to a third party logistics company, while ensuring that all delivery tasks are performed.
Customer orders have to be delivered as quickly as possible via direct delivery, and multiple orders cannot be combined within a single trip by an occasional driver, as in \citet{Archetti2016, ausseil2022supplier}.  

Our main contributions can be summarized as follows:
\begin{itemize}
\item We introduce the \probabbr, which considers (generic) compensation-dependent task acceptance probability functions of occasional drivers and allows to derive optimal compensations offered to them while minimizing the expected cost to the operator.
\item We introduce a MINLP formulation for the \probabbr\ (with generic acceptance probability functions) and show that it can be solved optimally via a two-stage approach that decomposes compensation and assignment decisions. We also show, that this two-stage approach which uses an exact linearization of our MINLP formulation allows to solve the \probabbr\ in polynomial time for generic probability functions under mild conditions that are detailed in \cref{sec:apuc}.
\item We study two practically relevant classes of acceptance probability functions (linear and logistic acceptance probability functions) and derive explicit formulas for optimal compensation values in these cases. 
\item We study a generalization of the \probabbr, show that it is NP-hard, and propose an approximate
linearization scheme for an appropriately extended MINLP formulation. 
\item We conduct an extensive computational study and sensitivity analysis on the main model parameters. Our results show that the compensation scheme proposed in this work consistently and significantly outperforms alternative established compensation schemes from the literature. 
We also demonstrate that this scheme ensures higher satisfaction of occasional drivers by providing them with more and more successful offers. Finally, the extremely low runtimes on large instances consisting of 100 tasks and up to 150 occasional drivers and the computational complexity of the proposed method indicate that our approach can be applied to larger problem instances that could arise in real-life applications.
\item We show how to use our two-stage approach as a heuristic in settings in which drivers and tasks appear dynamically over time. Our computational results indicate the usefulness of our two-stage method in dynamic settings and that the compensation scheme proposed in this work outperforms alternatives for this case as well.
\end{itemize}

This paper is organized as follows. \cref{sec:literature} summarizes the related literature. \cref{sec:apuc} formally introduces the problem studied in this paper, provides a MINLP formulation, and discusses the considered acceptance probability functions. \cref{sec:ilp} introduces the theoretical results that lead to the exact linearizations 
and polynomial-time solvability of \probabbr. \cref{sec:non-separable} focuses on the aforementioned generalization of the \probabbr, shows that this variant is NP-hard, and introduces an approximate linearization scheme for this more general case. \cref{sec:experimental-setup} details the setup of our computational study, whose results and findings are discussed in \cref{sec:results}. We conclude in \cref{sec:conclusion}, where we also provide possible future research directions. 
The Appendix contains proofs of all theoretical results, and additional computational results
\ifTR
.
\fi

\section{Literature review} \label{sec:literature}

In this section, we review articles that address uncertainty in driver behavior and/or compensation optimization subjects in the context of crowdsourced delivery. The reviewed papers consist of those identified by \citet{savelsbergh2022challenges}, as well as more recent papers addressing both subjects. 
\cref{table:literature} provides an overview of these articles and classifies them 
according to three criteria: service type
(following the classification of \citet{sampaio_chapter_2019}), 
compensation, and crowdshipper acceptance probability. An article is classified as driver dependent (or independent) if the attributes of occasional drivers are considered (or not considered) in the corresponding criterion.  
Integrated compensation indicates that
compensation decisions are simultaneously made with other decisions. Finally, works considering compensation dependent acceptance probability model acceptance probability of occasional drivers as a function of compensation.
In the following paragraph, we briefly summarize the articles that consider uncertainty in driver behavior without using acceptance probability functions. \cref{sec:compensation} discusses papers in which the compensation offered to drivers for performing tasks are treated as decisions that are independent of acceptance probabilities. Finally, \cref{sec:acceptance_uncertainty_and_compensation} provides details on articles that assume a relationship between compensation decisions and acceptance probabilities.

\begin{table}
	\centering
		\caption{Literature overview}
		\begin{tabular}{l|cc|ccc|ccc}
			\toprule
			&  \multicolumn{2}{c}{Service Type}   &                                                \multicolumn{3}{c}{\parbox{7em}{\centering Compensation}}                                                &                        \multicolumn{3}{c}{\parbox{7em}{\centering Acceptance Probability}}                         \\
			\cmidrule(lr){2-3} \cmidrule(lr){4-6} \cmidrule(lr){7-9}
			Paper & \rot{Multi-Depot} & \rot{In-Store} &  \rot{\parbox{7em}{\raggedright Driver Independent}} & \rot{\parbox{7em}{\raggedright Driver Dependent}}  & \rot{\parbox{7em}{\raggedright Integrated}} & \rot{\parbox{7em}{\raggedright Driver Independent}} & \rot{\parbox{7em}{\raggedright Driver Dependent}} & \rot{\parbox{7em}{\raggedright Compensation Dependent}} \\ \midrule
			\citet{kafle2017design}                                          &       \OK              &            &            &                  &                                                                          &                                                     &                                                   &                                                         \\
			\citet{gdowska2018stochastic}                                    &                    &      \OK       &            &                  &                                                                          &                         \OK                         &                                                   &                                                         \\
			\citet{allahviranloo2019dynamic}                                 &        \OK         &                &            &                  &                                                                          &                                                     &                                                   &                                                         \\
			\citet{mofidi2019beneficial}                                     &        \OK         &                &            &                  &                                                                          &                                                     &                                                   &                                                         \\
			\citet{dai2020workforce}                                         &         \OK            &            &            &                  &                                                                          &                                                     &                                                   &                                                         \\
			\citet{ausseil2022supplier}                                      &        \OK         &                &            &                  &                                                                          &                                                     &                                                   &                                                         \\
			\citet{behrendt2022prescriptive}                                 &        \OK         &                &            &                  &                                                                          &                                                     &                                                   &                                                         \\
			\citet{santini2022probabilistic}                                 &                    &      \OK       &            &                  &                                                                          &                         \OK                         &                                                   &                                                         \\
			\citet{dahle2019pickup}                                          &        \OK         &                &            &       \OK        &                                   \OK                                    &                                                     &                                                   &                                                         \\
			\citet{yildiz2019service}                                        &      \OK               &            &    \OK     &                  &                                                                          &                         \OK                         &                                                   &                                                         \\
			\citet{cao_last-mile_2020}                                       &         \OK            &            &            &       \OK        &                                                                          &                                                     &                                                   &                                                         \\
			\citet{le2021designing}                                          &        \OK         &                &            &       \OK        &                                   \OK                                    &                                                     &                                                   &                                                         \\
			\citet{Barbosa2022}                                              &                    &      \OK       &    \OK     &                  &                                                                          &                         \OK                         &                                                   &                           \OK                           \\
			\citet{hou_optimization_2022}                                    &                    &      \OK       &            &       \OK        &                                                                          &                                                     &                        \OK                        &                           \OK                           \\
			Our Work                                                         &        \OK         &      \OK       &            &       \OK        &                                   \OK                                    &                                                     &                        \OK                        &                           \OK                           \\ \bottomrule
		\end{tabular}
		\label{table:literature}
\end{table}

\citet{kafle2017design} and \citet{allahviranloo2019dynamic} focus on a bid submission setting in which a planner obtains bids that contain information about tasks and compensation amounts required by occasional drivers. \citet{mofidi2019beneficial} and \citet{ausseil2022supplier} consider approaches in which a planner offers menus of tasks to occasional drivers. \citet{behrendt2022prescriptive} study a scheduling problem that considers occasional drivers who are willing to commit to a time slot. \citet{dai2020workforce} focus on a workforce allocation problem that considers the allocation of occasional drivers to different restaurants in a meal delivery setting. Finally, \citet{gdowska2018stochastic} and \citet{santini2022probabilistic} study settings in which the acceptance probabilities of occasional drivers are considered but represented by fixed values.

\subsection{Independent compensation decisions} \label{sec:compensation}

\cite{dahle2019pickup} propose an extension of the pick-up and delivery vehicle routing problem with time windows in which the company owning the fleet can use occasional drivers to outsource some requests. 
They model the behavior of occasional drivers using personal threshold constraints, which indicate the minimum amount of compensation for which an occasional driver is willing to accept a task, and study the impact of the following three compensation schemes: fixed and equal compensation for each served request, compensation proportional to the cost of traveling from the pickup location to the delivery location, and compensation proportional to the cost of the detour taken by the occasional driver. 
Their results show that using occasional drivers can lead to savings of about 10-15\%, even when using a sub-optimal compensation scheme. They also show that these savings can be further increased by using more complex compensation schemes. 

\cite{yildiz2019service} model a crowdsourced meal delivery system as a multi-server queue with general arrival and service time distributions to investigate which restaurants should be included in the network, what payments should be offered to couriers, and whether full-time drivers should be used. 
They study the trade-off between profit and service quality and analyze the impact of unit revenue and vehicle speed. 
They derive optimal compensation amounts analytically, considering per-delivery and per-mile schemes, with respect to the optimal service area. They also investigate the impact of a given acceptance probability on the service area and compensation.

\citet{cao_last-mile_2020} consider a 
sequential packing problem to model task deliveries by occasional drivers and a professional fleet. In their setting, a distribution center receives delivery tasks at the beginning of a day, assigns bundles of tasks to occasional drivers during the day, and delivers the remaining tasks by the company-owned fleet at the end of the day. To model the arrival of occasional drivers, a marked Poisson process is proposed that specifies the time at which an occasional driver requests a task bundle and the bundle size. The arrival rate is assumed to be the same for each task bundle and is influenced by an incentive rate that allows an occasional driver to earn the same profit by delivering any bundle. They optimize the incentive rate that minimizes total delivery costs.

\cite{le2021designing} integrate task assignment, pricing, and compensation decisions. They assume that an occasional driver will accept a task if the compensation is above a given threshold and below an upper bound that the operator is willing to offer. These (distance-dependent) thresholds are derived from real-world surveys. The authors evaluate flat and individual compensation schemes under different levels of supply and demand. 

\subsection{Compensation dependent acceptance probability} \label{sec:acceptance_uncertainty_and_compensation}

\citet{Barbosa2022} extend the problem proposed by \citet{Archetti2016} by integrating compensation decisions. As one of the few papers to date, they combine pricing, matching, and routing aspects while considering the possibility that occasional drivers may refuse tasks.  
The acceptance probability is modeled as a function of the compensation offered. However, the authors assume that the probability function is identical 
for each pair of occasional driver and task. This implies that decisions about compensations are reduced to defining a single compensation that is constant for each occasional driver offered a task.
The considered problem is solved using a heuristic algorithm which builds upon the method proposed in \citet{gdowska2018stochastic}.

\citet{hou_optimization_2022} propose an optimization framework for a crowdsourced delivery service. They model the problem as a discrete event system, where an event is a task or a driver entering the system. A two-stage procedure is proposed to solve the problem. Tasks are assigned in the first stage, which focuses on minimizing the total detour of assigned tasks. 
In the second stage, optimal compensations are calculated that minimize the total delivery cost for the given task. 
The problem setting also considers the possibility that occasional drivers may refuse tasks, and the acceptance probabilities are determined by a binomial logit discrete choice model.

Studying the impact of compensation decisions on the acceptance behavior of occasional drivers has been identified as an important research direction. The first steps towards filling this research gap have been taken by \citet{Barbosa2022} and \citet{hou_optimization_2022}. \citet{Barbosa2022} do not consider individual properties of driver-task pairs and solve the resulting problem heuristically. Even though \citet{hou_optimization_2022} consider compensations for each assigned pair individually, their sequential algorithm computes optimal compensations once the assignments are fixed. As a result, their approach lacks an evaluation of all possible task assignment-compensation combinations. We aim to overcome these shortcomings by making assignment and compensation decisions simultaneously while considering compensation-dependent acceptance behaviors of occasional drivers, and by proposing an exact method for solving the resulting optimization problem. Additionally, most studies in the crowdsourced delivery literature consider either a multi-depot setting, which consists of tasks with different origins and destinations, or an in-store delivery setting in which crowdshippers are customers originating from a single store and make deliveries en-route. Our model is capable of solving problems in both environments.

\section{Problem definition\label{sec:apuc}} 

In the following, we formally introduce the \probname. 
The \probabbr aims to fulfill a set of online orders (\emph{tasks}) $I$ by offering them to occasional drivers or outsourcing them to a third part logistics company (3PL). 
The set of occasional drivers $J$ represents ordinary customers and other independent couriers who express willingness to perform tasks from $I$. These drivers are not fully employed and perform individual delivery tasks on their own time. 
A solution to the \probabbr allocates each task $i\in I$ either by offering the task to an occasional driver $j\in J$ or by 
outsourcing it to the 3PL who is assumed to perform the task at a costs $c_i\ge 0$.
In the remainder, we refer to the 3PL as professional drivers to differentiate them from the occasional drivers.

Occasional drivers are offered a compensation 
for performing a task, and they can choose whether to accept the task given their own utility with regard to compensation.
The compensation dependent \emph{acceptance probability} $P_{ij}(C_{ij})$ specifies the probability that an occasional driver $j\in J$ will accept an offer to perform a task $i\in I$ when offered a compensation of $C_{ij} \ge 0$. Tasks refused by an occasional driver must be performed by a professional driver, incurring a penalized cost $c_i'\ge c_i$ to the company. 

A feasible solution of the \probabbr can be described as a one-to-one assignment $A \subset I_\mathrm{o} \times J_\mathrm{o}$  where $I_\mathrm{o}\subseteq I$ is the set of tasks offered to a subset of the occasional drivers $J_\mathrm{o}\subseteq J$. Every other task $i \in I  \setminus I_\mathrm{o}$ is performed by a professional driver, while no tasks are offered to the occasional drivers belonging to $J \setminus J_\mathrm{o}$. The total expected cost of such a solution is represented by
\begin{equation}
\sum_{i\in I \setminus  I_\mathrm{o}} c_i + \sum_{ (i,j) \in A} \left(P_{ij}(C_{ij}) C_{ij} + (1-P_{ij}(C_{ij})) c_i' \right), \label{eq:obj}
\end{equation}
which sums up the total delivery cost to perform each task $i\in I \setminus I_\mathrm{o}$ by professional drivers and
the expected cost stemming from offering compensation $C_{ij}$ to driver $j$ to perform task $i$, for each $(i,j) \in A$.
The expected cost for each $(i,j) \in A$ is given by the compensation $C_{ij}$ times the probability $P_{ij}(C_{ij})$ that occasional driver $j$ will accept task $i$, plus 
the penalized cost
$c'_{i}$ to perform task $i$, times the probability $1 - P_{ij}(C_{ij})$ that occasional driver $j$ will refuse the task.
The objective of the \probabbr is to define the sets $I_\mathrm{o}$ and $J_\mathrm{o}$, their one-to-one assignment $A \subset I_\mathrm{o} \times J_\mathrm{o}$, and the compensation $C_{ij}$ for each $(i,j) \in A$ in such a way that the expected total cost is minimized.

\begin{figure}[!ht]
\begin{subfigure}[b]{.30\linewidth}
	\begin{tikzpicture}
		\node[anchor=south west,inner sep=0] (image) at (0,0) {\includegraphics[width=\textwidth]{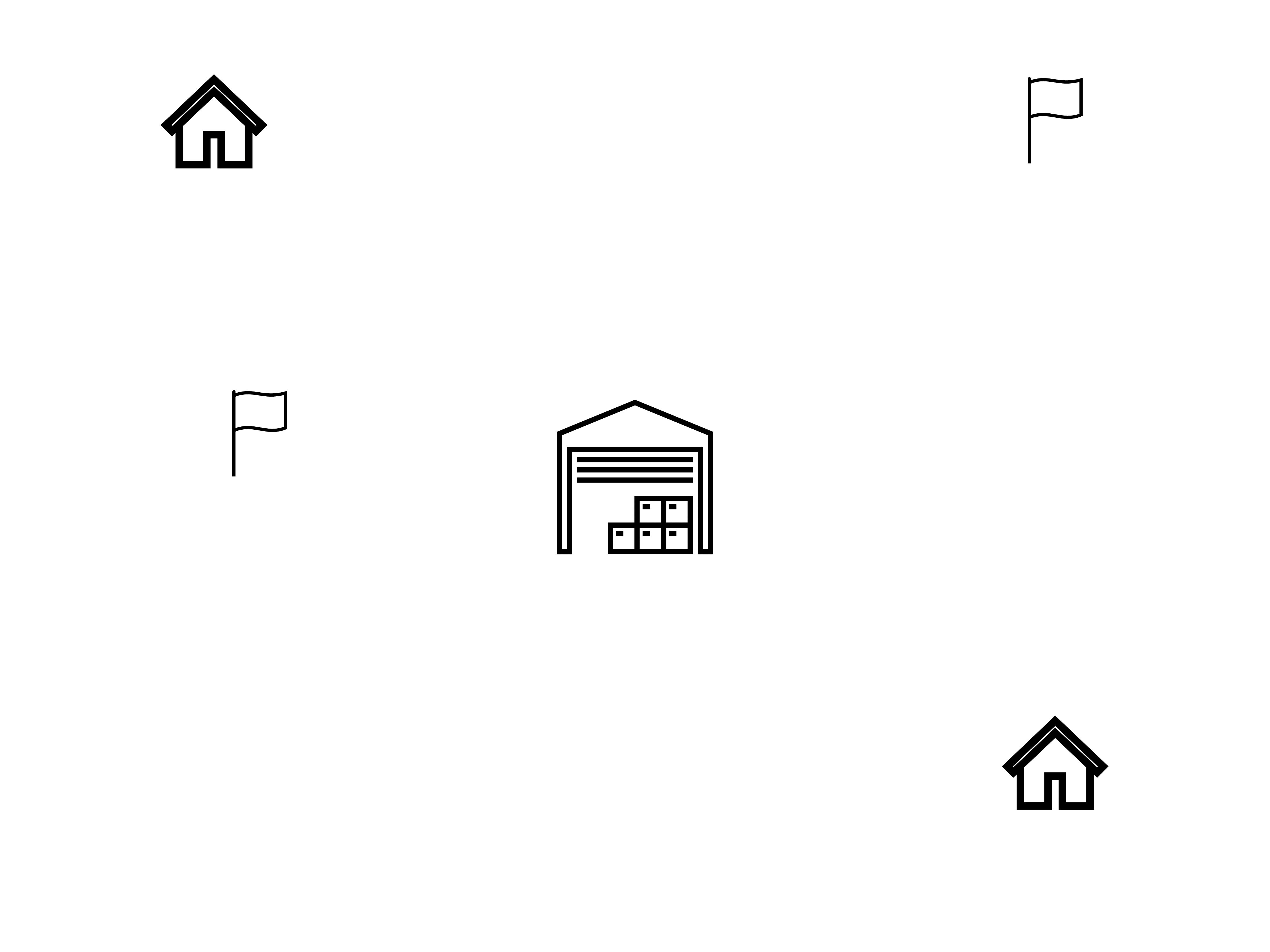}};
		\begin{scope}[x={(image.south east)},y={(image.north west)}]
			\node (task1) at (0.84,0.95) {\tiny Task 1};
			\draw[red, thick](0.84,0.9) ellipse (0.1 and 0.1);
			\node (company cost) at (0.84,0.75) {\tiny $c_1$};

			\node (task2) at (0.2,0.475) {\tiny Task 2};			
			\node (od1) at (0.175,0.95) {\tiny OD 1};
			\draw[rotate=7,red,  dashed](0.25,0.67) ellipse (.15 and 0.32);
			\node[blue] (comp) at (0.35,1) {\tiny $C_{21}$}; 
		
		\node (od2) at (0.835,0.28) {\tiny OD 2};
		\draw[thick, ->](0.57,0.42) -- (0.78,0.21);
		
	\end{scope}
\end{tikzpicture}
\caption{Instance and solution.} 
\label{fig:example_a}
\end{subfigure}
\hfill
\begin{subfigure}[b]{.30\linewidth}
\begin{tikzpicture}
\node[anchor=south west,inner sep=0] (image) at (0,0) {\includegraphics[width=\textwidth]{Figures/Figures2.png}};
\begin{scope}[x={(image.south east)},y={(image.north west)}]
	\node (task1) at (0.84,0.95) {\tiny Task 1};
	\draw[red, thick](0.84,0.9) ellipse (0.1 and 0.1);
	\node (company cost) at (0.84,0.75) {\tiny $c_1$};

	\node (task2) at (0.2,0.475) {\tiny Task 2};			
	\draw[blue,thick, ->](0.425,0.5) -- (0.24,0.55);
	\node (od1) at (0.175,0.95) {\tiny OD 1};
	\draw[rotate=7,red,  dashed](0.25,0.67) ellipse (.15 and 0.32);
	\draw[blue,thick, ->](0.2,0.6) -- (0.175,0.8);
	\node[blue] (comp) at (0.4,1) {\tiny $\checkmark$}; 
	\node[blue] (comp) at (0.4,.925) {\tiny $P(C_{21})$}; 
	\node[blue] (comp) at (0.4,.85) {\tiny $C_{21}$}; 

\node (od2) at (0.835,0.28) {\tiny OD 2};
\draw[thick, ->](0.57,0.42) -- (0.78,0.21);

\end{scope}
\end{tikzpicture}
\caption{OD~1 accepts task 2.} 
\label{fig:example_b}
\end{subfigure}
\hfill
\begin{subfigure}[b]{.30\linewidth}
\begin{tikzpicture}
\node[anchor=south west,inner sep=0] (image) at (0,0) {\includegraphics[width=\textwidth]{Figures/Figures2.png}};
\begin{scope}[x={(image.south east)},y={(image.north west)}]
\node (task1) at (0.84,0.95) {\tiny Task 1};
\draw[red, thick](0.84,0.9) ellipse (0.1 and 0.1);
\node (company cost) at (0.84,0.75) {\tiny $c_1$};

\node (task2) at (0.2,0.475) {\tiny Task 2};	
\draw[red, thick](0.2,0.525) ellipse (0.1 and 0.1);
\node (penalty cost) at (0.2,0.675) {\tiny $c'_2$};	
\node (od1) at (0.175,0.95) {\tiny OD 1};
\draw[blue,thick, ->](0.475,0.6) -- (0.21,0.85);
\draw[rotate=7,red,  dashed](0.25,0.67) ellipse (.15 and 0.32);
\node[blue] (comp) at (0.4,1) {\tiny \crossmark}; 
\node[blue] (comp) at (0.45,.925) {\tiny $1-P(C_{21})$}; 

\node (od2) at (0.835,0.28) {\tiny OD 2};
\draw[thick, ->](0.57,0.42) -- (0.78,0.21);

\end{scope}
\end{tikzpicture}
\caption{OD~1 refuses task 2.}
\label{fig:example_c}
\end{subfigure}
\caption{An instance and a solution of the \probabbr in which occasional driver (OD)~1 is offered a compensation of $C_{21}$ for task~2. Figures~\ref{fig:example_b} and \ref{fig:example_c} visualize the deliveries when OD~1 accepts or refuses task~2, respectively.} \label{fig:Example}
\end{figure}

\cref{fig:Example} provides an instance of the \probabbr and illustrates the allocation of tasks to occasional and professional drivers,
as well as the role of the acceptance probability.
Within this example, we consider a scenario with two tasks and two occasional drivers, i.e., $I = \{ 1, 2\}$ and $J = \{ 1, 2\}$. The figure indicate the store (in the middle of the figure) and the delivery location of the two tasks as well as the home location of the occasional drivers (to illustrate the potential detour of potential occasional drivers to complete the task). For simplicity, we assume that both occasional drivers are initially located at the central store.
\cref{fig:example_a} illustrates a solution for the instance in which task 1 is assigned to a professional driver (at cost $c_1$) and compensation $C_{21}$ is offered to occasional driver 1 for performing task 2, i.e., $I_{\mathrm{o}} = \{ 2\}$, $J_{\mathrm{o}} = \{ 1\}$, and $A = \{(2,1)\}$. Occasional driver 2 drives directly from the central depot to its final destination while occasional driver 1 must accept or refuse the task. With probability $P(C_{21})$, occasional driver 1 will accept the task (leading to total costs of $C_{21}$ + $c_1$) in which case they perform the task and then drive to their final destination, see \cref{fig:example_b}. With probability $1 - P(C_{21})$, they refuse the task and drive directly to their final destination. In this case, an additional professional driver must be used to perform task 1, and the total costs are equal to $c_1 + c_2'$, see \cref{fig:example_c}. Therefore, the expected total cost of the solution equals $P(C_{21})(C_{21} + c_1) + (1 - P(C_{21}))(c_1 + c'_2)$.

\medskip
\paragraph{Assumptions}
The above definition of the \probabbr assumes a static scenario with complete information, where the set of tasks $I$ and the set of available occasional drivers $J$ are known at the time of the assignment (e.g., by requiring the occasional drivers to sign up for a certain time slot in order to be offered tasks). While the availability of an occasional driver is certain, their acceptance behavior is probabilistic. 
Within our framework, we only 
assume that
the underlying compensation-dependent acceptance probability functions $P_{ij}: \R_+ \mapsto [0,1]$ are \emph{separable} for each $i\in I$ and $j\in J$,
i.e., that the probability of driver $j\in J$ for task $i\in I$ does not depend on other tasks $i'\in I\setminus \{i\}$ or compensations offered to other drivers $j'\in J\setminus \{j\}$, nor on their acceptance decisions. We argue that the latter assumption is realistic in the considered static setting where occasional drivers receive individual offers and thus have little opportunity to exchange information. 
Notice that information exchange between occasional drivers after their acceptance or refusal, which may influence their behavior, can be incorporated into the acceptance probabilities before each time slot in a dynamic setting while obeying the separability assumption.
Finally, we also assume that $P_{ij}(0) = 0$ for all $i\in I, j \in J$, i.e., no occasional driver works for free and, consequently, is never offered a task without compensation. In practice, we assume that historical information on driver decisions combined with characteristics of tasks and occasional drivers (e.g., their destinations or historical acceptance behavior), and external factors allow to estimate the probability function.

\subsection{Mixed integer nonlinear programming formulation\label{sec:mnlp}}

Next, we introduce an MINLP formulation \eqref{eq:minlp} for the \probabbr that uses three sets of decision variables. 
Variables $y_i\in \{0,1\}$ indicate whether task $i\in I$ is performed by a professional driver, and variables $x_{ij}\in \{0,1\}$ indicate whether task $i\in I$ is offered to occasional driver $j\in J$. Variables $C_{ij}\ge 0$ represent the compensation offered to occasional driver $j\in J$ for performing task $i\in I$.

\begin{subequations}\label{eq:minlp}
\begin{align}
\min \quad   & \sum_{i\in I} \left( c_i y_i +  \sum_{j\in J} \left( P_{ij}(C_{ij})C_{ij}  + (1-P_{ij}(C_{ij})) c'_i \right) x_{ij}   \right) \label{eq:minlp:obj} \\
\mbox{s.t.}\quad & \sum_{i\in I} x_{ij} \leq 1 & j\in J \label{eq:minlp:oc:assignment} \\
& \sum_{j\in J} x_{ij} + y_i = 1 & i\in I \label{eq:minlp:task:assignment} \\
& C_{ij} \le U_{ij} x_{ij} & i\in I,\ j\in J \label{eq:minlp:forcing} \\
& x_{ij} \in \{0,1\} & i\in I,\ j\in J \label{eq:minlp:x}  \\
& y_i \in \{0,1\} & i\in I \label{eq:minlp:y}  \\
& C_{ij} \ge 0 & i\in I,\ j\in J \label{eq:minlp:C}
\end{align}
\end{subequations}

The objective function~\eqref{eq:minlp:obj} minimizes the expected total cost. 
Constraints~\eqref{eq:minlp:oc:assignment} ensure that at most one task is offered to each occasional driver. Equations~\eqref{eq:minlp:task:assignment} ensure that each task is either offered to an occasional driver or performed by a professional driver. 
Inequalities~\eqref{eq:minlp:forcing} ensure that the compensation offered is not greater than a given upper bound $U_{ij}$ if task $j\in J$ is offered to occasional driver $i\in I$. Note that the cost $c_i$ of a professional driver for task $i\in I$ can always be used as such an upper bound, since it will never be optimal to offer compensation that exceeds this value. Furthermore, constraints~\eqref{eq:minlp:forcing} also force $C_{ij}$ to zero if task $i\in I$ is not offered to occasional driver $j\in J$. 
Finally, constraints \eqref{eq:minlp:x}--\eqref{eq:minlp:C} define the domains of the variables.
The difficulty of formulation~\eqref{eq:minlp} depends significantly on the considered probability function $P_{ij} : \R_+ \mapsto [0,1]$ of occasional drivers $j\in J$ for tasks $i\in I$. In \cref{sec:ilp}, we discuss a two-stage solution framework considering a general acceptance probability function. In the remainder of this section, we introduce two probability functions that can be used to model realistic acceptance behavior and discuss their use within the solution framework in \cref{sec:ilp}. 


\subsection{Linear acceptance probability\label{sec:mnlp:linear}}
One way to model the acceptance behavior of occasional drivers is to assume that their willingness to perform a task increases linearly with the compensation offered. \citet{campbell_incentive_2006} propose a linear approach to model the time slot selection behavior of customers in the context of attended home delivery services. We adopt their idea to model the acceptance behavior of occasional drivers. 
For each task-driver pair, we assume an initial acceptance probability that depends on the pair's characteristics. We also assume that this probability can be increased linearly by the amount of compensation offered, at a rate that also depends on the pair's characteristics. Besides task specific characteristics such as distance or weight, these can also include other factors such as weather and traffic levels. Such linear acceptance models allow avoiding the oversimplified assumption of a given (task-specific) threshold above which an occasional driver will always perform a task. A linear model also has the advantage of having a small number of parameters that need to be estimated. 
We consider a linear model in which we assume that occasional driver $j\in J$ accepts task $i\in I$ with a given \emph{base probability} $\alpha_{ij}$, $0\le \alpha_{ij}\le 1$, if the compensation offered is greater than zero. 
Parameter $\beta_{ij}> 0$ specifies the rate at which the acceptance probability of occasional driver $j\in J$ increases for task $i\in I$ for a given compensation $C_{ij}$. Hence, the probability function is formally defined as
\begin{equation} \label{eq:prob:linear}
P_{ij}(C_{ij}) = \begin{cases}
0 & \mbox{if $C_{ij}=0$} \\
\min \{\alpha_{ij} + \beta_{ij} C_{ij}, 1 \} & \mbox{ otherwise }			
\end{cases}, \quad i\in I,\ j\in J.
\end{equation}

Considering function~\eqref{eq:prob:linear} in formulation~\eqref{eq:minlp}, we first observe that there always exists an optimal solution in which $C_{ij}>0$ if task $i\in I$ is offered to occasional driver $j\in J$, i.e., if $x_{ij}=1$. Since the acceptance probability for a compensation of zero is equal to zero, the alternative solution in which task $j$ is performed by a professional
driver is always at least as good (recall that we assumed that $c_i'\ge c_i$). Similarly, it is never optimal to offer a compensation that exceeds the cost of a professional driver for the same task, or that exceeds the minimum compensation that leads to an acceptance probability equal to one. Thus, $\min \{c_i, \frac{1 - \alpha_{ij}}{\beta_{ij}}\}$ is an upper bound on the compensation offered to occasional driver $j\in J$ for task $i\in I$ that can be used to (possibly) strengthen the upper bounds $U_{ij}$ in constraints~\eqref{eq:minlp:forcing}.
Consequently, the variant of formulation~\eqref{eq:minlp} for the linear acceptance probability case is obtained by replacing the objective function~\eqref{eq:minlp:obj} by
\begin{align}
\min  \quad  &  \sum_{i\in I} \left( c_i y_i +  \sum_{j \in J} [(\alpha_{ij}+\beta_{ij} C_{ij}) C_{ij}  + (1-\alpha_{ij}-\beta_{ij} C_{ij}) c'_i] x_{ij}   \right). \label{eq:mnlp:obj:linear-prob}
\end{align} 

\subsection{Logistic acceptance probability\label{sec:mnlp:logistic}}

An alternative to modeling the acceptance behavior of occasional drivers that has been suggested in the literature is the use of logistic regression models. \citet{devari_crowdsourcing_2017} propose a logistic regression model to model the acceptance behavior of occasional drivers for performing a delivery for their friends. Similarly, \citet{Le2019} propose a binary-logit model where the willingness to work as an occasional driver is considered a dependent variable.
Logistic regression is particularly appealing when historical or experimental data are available on task and driver attributes, external factors, the compensation offered for tasks and the resulting acceptance decisions. Such data can be used to fit a logistic regression model that predicts the acceptance probabilities of occasional drivers. 
Using parameters $\gamma_{ij}$, which subsumes the impact of all predictors other than the compensation, and $\delta_{ij}>0$, logistic acceptance probabilities are modeled using functions 

\begin{equation} \label{eq:prob:logistic}
P_{ij}(C_{ij}) = \begin{cases}
0 & \mbox{if $C_{ij} = 0$} \\
\frac{1}{1 + e^{ - (\gamma_{ij} + \delta_{ij} C_{ij})}}	 & \mbox{ otherwise}
\end{cases}, \quad i\in I,\ i\in J.
\end{equation}

A nonlinear model for the \probabbr when considering logistic acceptance probabilities is obtained from \eqref{eq:minlp} if the objective function \eqref{eq:minlp:obj} is replaced by 
\begin{align}
\min  \quad  &  \sum_{i\in I} \left( c_i y_i +  \sum_{j \in J} [(\frac{1}{1 + e^{ - (\gamma_{ij} + \delta_{ij} C_{ij})}}) C_{ij}  + (1-\frac{1}{1 + e^{ - (\gamma_{ij} + \delta_{ij} C_{ij})}}) c'_i] x_{ij}   \right). \label{eq:mnlp:obj:logistic-prob}
\end{align}

\subsection{Learning acceptance probabilities from historical data}

Crowdsourced delivery platforms use online tools to process delivery orders and allocate them to occasional drivers. This enables them to collect a wide range of data about their operations. In the following, we discuss how historical data about offers made to occasional drivers can be used to learn acceptance probabilities that can be used within our framework. For this, the data collected for previous offers must contain information about 
\begin{inparaenum}[(i)]
\item the task,
\item the driver, 
\item the compensation offered to the driver, and
\item whether the driver accepted the offer
\end{inparaenum}
as illustrated in \cref{fig:snapshot}.

The required information associated with tasks and drivers depends on the operational setting and considered probability function. 
For the in-store setting described in the previous sections, it is sufficient to keep track of task delivery points and driver destinations to calculate the detour, since both, task shipment points of the tasks as well as driver initial locations are assumed to be the store. 
In a setting where tasks have distinct shipment points and drivers have distinct initial locations, keeping track of such additional information is required to calculate the detour. 
As described below, additional information about the environment must be included in the database in case their impact should be incorporated in the acceptance probability functions.

\begin{table}[H]
	\caption{A snapshot of representative operational data} 
	\resizebox{\textwidth}{!}{%
		\begin{tabular}{lllllllll}
			\toprule
			\multicolumn{3}{c}{Task   information (i)} & \multicolumn{3}{c}{Driver information (ii)} & \multicolumn{3}{c}{Offer information}           \\ \cmidrule(lr){1-3} \cmidrule(lr){4-6} \cmidrule(lr){7-9}
			ID    & Shipment point  & Delivery point & ID      & Location    & Destination    & Compensation (iii) & Detour & Response (iv) \\ \midrule
			T001  & P001             & P003              & D001 & P005    & P006                         & 3.7 EUR      & 5.3 KM         & Reject          \\
			T002  & P002             & P004              & D002 & P002    & P007                         & 5.9 EUR      & 1.2 KM         & Accept          \\
			T003  & P001             & P004              & D003 & P001    & P008                         & 4.1 EUR      & 3.8 KM         & Reject          \\
			\dots  & \dots             & \dots             & \dots    & \dots           & \dots              & \dots      & \dots         & \dots          \\
			\bottomrule
		\end{tabular}%
	}
	\label{fig:snapshot}
\end{table}

A logistic (or linear) acceptance probability model can be developed using such historical data to estimate driver responses based on task and driver information, as well as offered compensation as independent variables. In such a model, the impact of compensation is captured by the $\delta$ (or $\beta$) parameter, while the combined influence of other factors is represented by the $\gamma$ (or $\alpha$) parameter in \cref{eq:prob:logistic} (or \cref{eq:prob:linear}). 
This approach also allows to easily incorporate additional factors related to drivers, orders, or the environment using the $\gamma$ ($\alpha$) parameter. 
Considering additional factors is especially important for revealing the indirect relationships that impact driver behavior. While the detour could be a significant factor that impact drivers' decisions during non-rush hours, the expected delivery time could, in contrast, be more significant during rush hours. Similarly, weather conditions can impact driver response depending on the vehicle used by the driver. A driver using bikes to make deliveries could be less willing to make extra deliveries on a rainy day. Hence, integrating such factors can improve predicting drivers' responses. In practice, such parameters can be easily integrated into the proposed methodology without changing the formulations.

Fitting accurate acceptance probability models at the individual level may be challenging due to potentially insufficient data for each individual. Instead, one can identify behavioral clusters among drivers and fit acceptance probability models that describe the behavior of each cluster. \citet{mohri_modeling_2024} demonstrates such an approach by employing a latent class model. Alternatively, one can use the complete dataset to fit an acceptance probability model that captures the effects of the observable factors by including them as variables in the model. Such an approach is demonstrated by \citet{Le2019} to evaluate the effects of various socioeconomic factors on drivers' decisions to participate in crowdsourced delivery.

\section{Mixed integer linear programming reformulations\label{sec:ilp}} 

A major advantage of the \probabbr over previously introduced models is its integration of assignment and compensation decisions, as well as its consideration of the uncertain acceptance behavior of occasional drivers. 
However, these aspects pose a significant challenge to the design of well-performing solution methods that derive optimal assignment and compensation decisions.
In this section, we show how to overcome this challenge by exploiting some of the assumptions made in the \probabbr. We will show that, under certain conditions, the \probabbr can be solved optimally using a two-phase approach that first identifies optimal compensation values that are used as input to the second phase, in which assignment decisions are made. We also provide explicit formulas for optimal compensation values in the case of linear and logistic acceptance probabilities. These results allow us to reformulate MINLP~\eqref{eq:minlp} and its variants for the latter two acceptance probability functions as mixed-integer linear programs (MILPs). They also imply that the \probabbr can be solved in polynomial time whenever the first-phase problem of identifying optimal compensation values can be solved in polynomial time.

The main property required for the results introduced below is that the compensation and acceptance decisions are \emph{separable}. This property holds for the \probabbr since the acceptance probability functions are assumed to be independent of each other (for each pair of task and occasional driver) and since the considered setting does not include constraints that limit the operator's choice, such as cardinality or budget limits on the tasks offered to occasional drivers.
\cref{th:compensation-values} reveals that the separability of compensation and acceptance decisions allows for the identification of optimal compensation values independent of the task allocation decisions. The latter is made more explicit in \cref{corr:compensations:optimal}.

\begin{theorem}\label{th:compensation-values}
	Consider an arbitrary instance of the \probabbr. There exists an optimal solution for this instance in which the (optimal) compensation values $C_{ij}^*$ are equal to
	\begin{equation*}
		C^*_{ij} = \argmin_{C_{ij}\ge 0} ~ P_{ij}(C_{ij})(C_{ij}-c'_i)
	\end{equation*}
	if task $i\in I$ is offered to occasional driver $j\in J$. Furthermore, compensations $C_{ij}^*$ for tasks $i\in I$ that are not offered to occasional driver $j\in J$ can be set to an arbitrary (non-negative) value.
\end{theorem}

\begin{proof}
	Provided in \ref{appendix:proofs}.
\end{proof}

%

\begin{corollary}\label{corr:compensations:optimal}
	There exists an optimal solution for an arbitrary instance of the \probabbr that can be identified by solving formulation \eqref{eq:minlp} while setting compensation values equal to
	\begin{equation}\label{eq:compensations:optimal}
		C_{ij}^* = \argmin_{C_{ij}\ge 0} P_{ij}(C_{ij}) (C_{ij} - c'_i)
	\end{equation}
	for each task $i\in I$ and every occasional driver $j\in J$.
\end{corollary}

A major consequence of \cref{corr:compensations:optimal} is that optimal compensation values, which can be calculated according to Equation~\eqref{eq:compensations:optimal}, are independent of the other decisions, i.e., which tasks to assign to professional drivers and which tasks to offer to which occasional driver. Therefore, we can first calculate the compensation values and then identify optimal assignment decisions using the MILP formulation~\eqref{eq:milp} in which cost parameters 
${\bar{c}}^*_{ij}$
are set to
\begin{equation}
	{\bar{c}}^*_{ij} = P_{ij}(C_{ij}^*) C_{ij}^* + (1 - P_{ij}(C^*_{ij})) c'_i  \label{eq:expected_cost}
\end{equation}
for each $i\in I$ and $j\in J$. 
As in formulation~\eqref{eq:minlp}, variables $y_i\in \{0,1\}$ indicate whether task $i\in I$ is performed by a professional driver and variables $x_{ij}\in \{0,1\}$ indicate whether it is offered to occasional driver $j\in J$. 
\begin{subequations}\label{eq:milp}
	\begin{align}
		\min \quad   & \sum_{i\in I} c_i y_i +  \sum_{i\in I} \sum_{j\in J} {\bar{c}}^*_{ij} x_{ij} \label{eq:milp:obj} \\
		\mbox{s.t.}\quad & \sum_{i\in I} x_{ij} \leq 1 & j\in J \label{eq:milp:oc:assignment} \\
		& \sum_{j\in J} x_{ij} + y_i = 1 & i\in I \label{eq:milp:task:assignment} \\
		& x_{ij} \in \{0,1\} & i\in I,\ j\in J \label{eq:milp:x}  \\
		& y_i \in \{0,1\} & i\in I \label{eq:milp:y}  
	\end{align}
\end{subequations}

Constraints~\eqref{eq:milp:oc:assignment} and \eqref{eq:milp:task:assignment} ensure that each occasional driver is offered at most one task, and that each task is either offered to an occasional driver or assigned to a professional driver. Overall, it is easy to see that formulation~\eqref{eq:milp} is equivalent to a standard assignment problem (with special ``assignment'' variables $y_i$ for each $i\in I$). Consequently, \cref{corr:polytime} follows, since formulation~\eqref{eq:milp} can be solved in polynomial time due to its totally unimodular constraint matrix.

\begin{corollary}\label{corr:polytime}
	The \probabbr can be solved in polynomial time if the compensation and acceptance decisions are separable and optimal compensation values $C^*_{ij}$ for tasks $i\in I$ and occasional drivers $j\in J$ can be identified in polynomial time.
\end{corollary}

\cref{prop:compenation-values:linear,prop:compenation-values:logistic} provide explicit formulas for optimal compensations in instances of the \probabbr if the acceptance behavior of occasional drivers is modeled using linear or logistic acceptance probability functions, respectively, as introduced in \cref{sec:apuc}. Thus, such instances satisfy the conditions of \cref{corr:polytime} and can be solved in polynomial time.

\begin{proposition}\label{prop:compenation-values:linear}
	Consider an instance of the \probabbr in which the acceptance behavior of occasional drivers is modeled using the linear acceptance probability function
	\begin{equation}\label{eq:prob:linear:proof} 
		P_{ij}(C_{ij}) = \begin{cases}
			0 & \mbox{if $C_{ij} =0$} \\
			\min \{\alpha_{ij} + \beta_{ij} C_{ij}, 1 \} & \mbox{ otherwise }			
		\end{cases}, \quad i\in I,\ j\in J.
	\end{equation}
	There exists an optimal solution for this instance with compensation values 
	$C_{ij}^*=\frac{c'_i}{2} - \frac{\alpha_{ij}}{2 \beta_{ij}}$
	for each task $i\in I$ and occasional driver $j\in J$.
\end{proposition}

\begin{proof}
	Provided in \ref{appendix:proofs}.
\end{proof}


\begin{proposition}\label{prop:compenation-values:logistic}
	Consider an instance of the \probabbr in which the acceptance behavior of occasional drivers is modeled using the logistic acceptance probability function 
	\begin{equation}\label{eq:prob:logistic:proof}
		P_{ij}(C_{ij}) = \begin{cases}
			0 & \mbox{if $C_{ij} = 0$} \\
			\frac{1}{1 + e^{ - (\gamma_{ij} + \delta_{ij} C_{ij})}}	 & \mbox{ otherwise}
		\end{cases}, \quad i\in I,\ i\in J.
	\end{equation}
	There exists an optimal solution to that instance with compensation values 
	$C_{ij}^* = - \frac{W(e^{\gamma_{ij} + \delta_{ij} c'_i - 1})-\delta_{ij} c'_i + 1}{\delta_{ij}}$
	for each task $i\in I$ and occasional driver $j\in J$, in which $W(\cdot)$ is the Lambert $W$ function. 
\end{proposition}

\begin{proof}
	Provided in \ref{appendix:proofs}.
\end{proof}

\section{The \probabbr without separability\label{sec:non-separable}}

In this section, we discuss a generalization of the \probabbr in which the decisions about which tasks to offer to occasional drivers and the corresponding compensation are not independent of each other. This occurs, for example, when the operator has to respect (strategic) considerations that limit the number of tasks offered to occasional drivers or the budget available for such offers. In this section, we assume that the relevant limitations can be modeled as a set of $L$ linear constraints involving assignment and compensation decisions. 
Using notation $a_{ij}^\ell$ and $b_{ij}^\ell$ for each $i\in I$, $j\in J$, and $\ell \in \{1, \dots, L\}$ to denote the coefficients associated to these two sets of variables and $B^\ell$ for the corresponding (resource) limits, the considered set of constraints is written as 
\begin{align}
	& \sum_{i\in I} \sum_{j\in J} ( a_{ij}^\ell x_{ij} + b_{ij}^\ell C_{ij}) \le B^\ell, & \ell\in \{1, \dots, L\}. \label{eq:non-separable}
\end{align}

A cardinality constraint on the total number of tasks offered to occasional drivers can, e.g., be realized by setting $a_{ij}^\ell=1$ and $b_{ij}^\ell=0$ for all $i\in I$ and $j\in J$ while $a_{ij}^\ell=0$ and $b_{ij}^\ell=1$ holds for a constraint limiting the overall budget offered to occasional drivers. 

\cref{prop:np} shows that the inclusion of \emph{non-separability constraints} \eqref{eq:non-separable} implies NP-hardness of the resulting variant of the \probabbr, which we call the \emph{\probname without Separability (\probnonsepabbr)}. 
	
	\begin{proposition}\label{prop:np}
		The \probnonsepabbr is strongly NP-hard.
	\end{proposition}

\begin{proof}
	Provided in \ref{appendix:proofs}.
\end{proof}

While considering the above set of constraints related to operator choices, we still assume that the acceptance decisions of occasional drivers depend only on the task and compensation offered, i.e., the acceptance probability functions $P_{ij}(C_{ij})$ are still applicable. Thus, while the results concerning optimal compensation decisions from \cref{sec:ilp} and the two-phase solution approach are no longer applicable, the objective function~\eqref{eq:minlp:obj} can be decomposed for each $i\in I$ and $j\in J$, which is exploited in \cref{prop:piecewise}.

\begin{proposition}\label{prop:piecewise}
	The nonlinear objective function~\eqref{eq:minlp:obj} in formulation~\eqref{eq:minlp} can be replaced by 
	\begin{align}
		\		& \sum_{i\in I} c_i y_i + \sum_{i\in I} \sum_{j\in J} (f_{ij}(C_{ij}) + g_{ij}(x_{ij})). \label{eq:minlp:obj:simplified}
	\end{align}
	Thereby, $f_{ij}(C_{ij})$ is a nonlinear function depending on $C_{ij}$, and $g_{ij}(x_{ij})$ is a linear function in $x_{ij}$. In the case of generic acceptance probability functions $P_{ij}(C_{ij})$, we have $f_{ij}(C_{ij})=P_{ij}(C_{ij})(C_{ij}-c_i')$ and $g_{ij}(x_{ij})=c_i' x_{ij}$.
\end{proposition}

\begin{proof}
	Provided in \ref{appendix:proofs}.
\end{proof}


\cref{prop:piecewise} implies that a piecewise linear approximation of the nonlinear objective function~\eqref{eq:minlp:obj} can be derived using standard techniques as described, e.g., in \citet{nemhauser_integer_1988}. For each task $i\in I$ and every driver $j\in J$, we approximate the non-linear function $f_{ij}(C_{ij})$ for $C_{ij}\in [0,U_{ij}]$ by a piecewise linear function defined using a discrete set of $K$ possible compensation values $u_{ij}^k$, $k=1, \dots, K$. Formally, the approximation is specified by the points $(u_{ij}^k, f_{ij}(u_{ij}^k))$, $k=1, \dots, K$, where $u_{ij}^1=0$, $u_{ij}^K=U_{ij}$, and $u_{ij}^{k-1}<u_{ij}^k$, $k\in \{2, \dots, K\}$.
For all $i\in I, j\in J$ and $k\in \{1, \dots, K\}$ we further define non-negative auxiliary variables $\lambda_{ij}^k\ge 0$ as the weights for the convex combination of the values of $f_{ij}({C}_{ij})$. Therefore we replace $f_{ij}({C}_{ij})$ in \eqref{eq:minlp:obj:simplified} by the convex combination $\sum_{k=1}^K \lambda_{ij}^k f_{ij}^k(u_{ij}^k)$ such that $\sum_{k=1}^K \lambda_{ij}^k=1$. Since the nonlinear functions $f_{ij}(C_{ij})$ are not necessarily convex, we need to additionally ensure that at most two of these variables are positive and that two such variables can only be both positive if they are consecutive, i.e., if $\lambda_{ij}^k>0$ and $\lambda_{ij}^r> 0$, then $r\in \{k-1,k+1\}$. To ensure these additional conditions, 
formulation~\eqref{eq:milp:approx} also uses binary variables $v_{ij}^k\in \{0,1\}$ for each $i\in I$, $j\in J$, and $k\in \{1, \dots, K-1\}$. Here, $v_{ij}^k=1$ indicates that $u_{ij}^k \le C_{ij} \le u_{ij}^{k+1}$. 
\begin{subequations}\label{eq:milp:approx}
	\begin{align}
		\min\quad & \sum_{i\in I} c_i y_i + \sum_{i\in I} \sum_{j\in J} \left( \sum_{k=1}^K  f_{ij}(u_{ij}^k) {\lambda_{ij}^k} + g_{ij}(x_{ij})\right) \label{eq:milp:approx-obj} \\	
		\mbox{s.t.}\quad  & \eqref{eq:minlp:oc:assignment}-\eqref{eq:minlp:y}, \eqref{eq:non-separable} \nonumber \\
		&\sum_{k=1}^K u_{ij}^k {\lambda_{ij}^k} \le U_{ij} x_{ij} & i\in I,\ j\in J \label{eq:milp:approx:U}\\
		& {\lambda_{ij}^1} \le v_{ij}^1 & i\in I,\ j\in J \label{eq:milp:approx:wv:1}\\
		& {\lambda_{ij}^k} \le v_{ij}^{k-1} + v_{ij}^k & i\in I, j\in J,\ k\in \{2, \dots, K-1\} \label{eq:milp:approx:wv}\\
		& {\lambda_{ij}^K} \le v_{ij}^{K-1} & i\in I,\ j\in J \label{eq:milp:approx:wv:K}\\		
		& \sum_{k=1}^{K-1} v_{ij}^k = 1 & i\in I,\ j\in J \label{eq:milp:approx:v:convex}\\
		& \sum_{k=1}^K {\lambda_{ij}^k} = 1 & i\in I,\ j\in J \label{eq:milp:approx:w:convex}\\
		& v_{ij}^k\in \{0,1\} & i\in I,\ j\in J,\ k\in \{1, \dots, K-1\}
		\label{eq:milp:approx:v} \\
		& {\lambda}_{ij}^k \ge 0 & i\in I,\ j\in J,\ k\in \{ 1,\dots, K \} \label{eq:milp:approx:w}\
	\end{align}
\end{subequations}

The objective function~\eqref{eq:milp:approx-obj} is a piecewise linear approximation of the one introduced in \cref{prop:piecewise} using the concepts described above.
Recall that $g_{ij}(x_{ij})$ is a linear function. As described in \cref{sec:mnlp}, constraints \eqref{eq:minlp:oc:assignment}-\eqref{eq:minlp:y} ensure that at most one task is offered to an occasional driver and that each task is either performed by a professional driver, or offered to an occasional driver. They also define the domains of variables $x_{ij}$ and $y_i$ for tasks $i\in I$ and occasional drivers $j\in J$. Constraints~\eqref{eq:non-separable} have been introduced at the beginning of this section while inequalities~\eqref{eq:milp:approx:U} are constraints $C_{ij}\le U_{ij} x_{ij}$ rewritten using the identity $\sum_{k=1}^K u_{ij}^k {\lambda}_{ij}^k = C_{ij}$. 
Finally, \eqref{eq:milp:approx:wv:1}--\eqref{eq:milp:approx:w} are standard constraints used to model piecewise linear approximations that ensure the conditions detailed above,
see, e.g., \citet{nemhauser_integer_1988} for further details. 

\medskip

\cref{prop:piecewise:linear,prop:piecewise:logistic} detail how formulation~\eqref{eq:milp:approx} can be modified when the acceptance behavior is modeled using a linear and logistic acceptance probability function, respectively.

\begin{proposition}\label{prop:piecewise:linear}
	Consider an arbitrary instance of the \probnonsepabbr in which the acceptance behavior of occasional drivers is modeled using the linear probability function~\eqref{eq:prob:linear}. Then, \cref{prop:piecewise} applies for $f_{ij}(C_{ij})=\beta_{ij} C_{ij}^2 + (\alpha_{ij} - \beta_{ij}c_i') C_{ij}$ and $g_{ij}(x_{ij}) = c_i'(1 - \alpha_{ij}) x_{ij}$.
	Furthermore, variables $v_{ij}^k$, $i\in I$, $j\in J$, $k\in \{1, \dots, K-1\}$ and constraints \eqref{eq:milp:approx:wv:1}--\eqref{eq:milp:approx:v:convex} involving them are redundant in formulation~\eqref{eq:milp:approx} and can therefore be removed.
\end{proposition}

\begin{proof}
	Provided in \ref{appendix:proofs}.
\end{proof}

%

\begin{proposition}\label{prop:piecewise:logistic}
	Consider an arbitrary instance of the \probnonsepabbr in which the acceptance behavior of occasional drivers is modeled using the logistic probability function~\eqref{eq:prob:logistic}. Then, \cref{prop:piecewise} applies for $g_{ij}(x_{ij}) =  c_i' x_{ij}$ and $f_{ij}(C_{ij})=\begin{cases}
		0 & \mbox{ if $C_{ij} = 0$} \\
		\frac{C_{ij} - c_i'}{1 + e^{-\gamma_{ij} - \delta_{ij} C_{ij}}} & \mbox{ otherwise}
	\end{cases}$. 
\end{proposition}

\begin{proof}
	Provided in \ref{appendix:proofs}.
\end{proof}


\section{Experimental setup\label{sec:experimental-setup}}
In this section, we describe our benchmark instances, provide further details on the parameters of the considered acceptance probability functions, and introduce alternative compensation models inspired by the literature that we use to evaluate our approach. 

\subsection{Instances \label{sec:instances}} 
As a relatively new area of research, the field of crowdsourced delivery still lacks established benchmark libraries. Most existing studies use instances from the VRP literature (e.g., \citet{Archetti2016,Barbosa2022}) or randomly generated synthetic instances (e.g., \citet{arslan_crowdsourced_2019,dayarian_crowdshipping_2020}). 
Since the instances used in these works do not consider acceptance probabilities, we generate a new set of synthetic instances. 
To this end, we simulate an in-store delivery setting, where occasional drivers are assumed to be regular in-store customers who are willing to deliver a task en route. We generate destinations for occasional drivers and tasks uniformly at random in a $200\times 200$ plane. The coordinates of these locations are rounded to two decimal places. The store is assumed to be located in the center of the plane and coincides with the initial locations of all occasional drivers. We assume that the 
professional delivery
cost $c_i$ for each task $i\in I$ is equal to the Euclidean distance between the task delivery point and the depot. 
Note that, since each task is identified with its origin and destination, the costs for professional and occasional drivers can be incorporated in the respective parameters without explicitly including these locations in the model. Therefore, our model generalizes problems consisting of tasks with different origins and destinations such as meal delivery, as well as tasks with single origin and different destinations such as in-store customer settings.
Each instance is characterized by a combination of the following three parameters:
\begin{enumerate}
	\item The number of occasional drivers $\odno \in \{50, 75, \dots, 150\}$.
	\item A penalty parameter $\pnlty \in \{0, 0.05, \dots, 0.25\}$ that defines the relative increase in 
	professional delivery
	costs if an occasional driver rejects an offer, i.e., $c_i'=(1+\rho) c_i$, $\forall i\in I$.
	\item A parameter $\utlty \in \{0, 0.1, \dots, 1\}$ that affects the acceptance probabilities of occasional drivers and, in particular, their sensitivity towards making detours when performing deliveries. In the linear case, this parameter is used to define the base probability $\alpha_{ij}$ for each $i\in I$ and $j\in J$ via a weighted distance utility $\nicefrac{\utlty d_{j}}{(d_{i} + d_{ij})}$, i.e., $\alpha_{ij}=\nicefrac{\utlty d_{j}}{(d_{i} + d_{ij})}$. Here, $d_k$ is the Euclidean distance between the store and the destination of $k\in I\cup J$, while $d_{ij}$ is the Euclidean distance between the destinations of task $i\in I$ and driver $j\in J$. 	
\end{enumerate}
For each configuration, five instances with 100 tasks are generated using different random seeds, so our instance library consists of $1\,650$ instances\footnote{The complete set of instances will be made available after publication.}.

\subsection{Acceptance probability functions} 

As explained above, for the \emph{linear model} we set the base probability $\alpha_{ij}$ of probability function~\eqref{eq:prob:linear} equal to the weighted distance utility for task $i\in I$ and driver $j\in J$. Thus, the base probability decreases with increasing detour of occasional drivers (for constant $\utlty$) or decreasing $\utlty$ (for constant detour). The rate of increase $\beta_{ij}$ is determined by dividing a random value from the uniform distribution $U(0.5,2)$ by
the detour for driver $j\in J$ when delivering task $i\in I$,
i.e. $\beta_{ij}=\nicefrac{X_j}{(d_i+d_{ij}-d_j)} : X_j \sim U(0.5,2)$. 

\smallskip
To simulate the \emph{logistic acceptance probabilities}, we first generate an artificial historical dataset containing $100\,000$ data points with simulated decisions of potential occasional drivers. Each data point consists of locations of a driver-task pair, $\alpha_{ij}$ and $\beta_{ij}$ values calculated as in the linear case, a compensation value drawn randomly from the uniform distribution over $[0,100\sqrt{2}]$ (the upper bound of the interval corresponds to the maximal 
professional delivery
cost that incurs for serving a customer request), and an acceptance decision generated using a Bernoulli trial based on the resulting (linear) acceptance probability.
Afterwards, this data set is used to train a logistic regression model in order to obtain parameters $\gamma_{ij}$ and $\delta_{ij}$ of the logistic acceptance probability function~\eqref{eq:prob:logistic}. The dependent variable considered in the logistic regression model is the acceptance decision, and the independent variables used to estimate it are the Euclidean distances between the store and the destinations of task and driver, the detour for delivering the task, the compensation, and the driver sensitivity. The driver sensitivity corresponds to $\beta_{ij}$ values used in the linear model.

\subsection{Benchmark compensation models \label{sec:benchmark_schemes}} 
We assess the potential advantages of the individualized compensation scheme proposed in this paper to the established detour-based, distance-based and flat compensation schemes introduced in the literature. The three schemes are formally defined as follows:

\begin{itemize}
\item In the \textbf{detour-based scheme}, the offered compensation is proportional to the detour of driver $j\in J$ when delivering task $i\in I$, i.e., $C_{ij}=p_\mathrm{detour} \cdot (d_i + d_{ij} - d_{j})$.
\item The \textbf{distance-based scheme} compensates proportional to the distance between the central store and the destination of task $i\in I$, i.e., $C_{ij} = p_\mathrm{distance} \cdot d_i$ holds for all $j\in J$.
\item The \textbf{flat scheme} offers a constant compensation for all tasks and drivers, i.e., $C_{ij} = p_\mathrm{flat}$ for all $i\in I$ and $j\in J$.
\end{itemize}

Each of the three schemes can be fully described by a single parameter (i.e., $p_\mathrm{detour}$, $p_\mathrm{distance}$, $p_\mathrm{flat}$). Naturally, the quality of the compensation schemes may be very sensitive to the values of these parameters. In order to assess the advantages of individual compensations over these schemes in a fair way, we aim to find (close-to) optimal parameter values for each of the three benchmark schemes. We use this value to replace the compensation in Equation~\eqref{eq:compensations:optimal} in the first phase of our two-phase solution approach, and then determine optimal assignments under the corresponding scheme. 

To determine the best compensation parameter for any given instance, we perform a heuristic search procedure in the interval $[0,p_\mathrm{max}]$ where $p_\mathrm{max}$ is the maximum value to which any optimal compensation $C_{ij}^*$ from the individualized scheme (identified using equation~\eqref{eq:compensations:optimal}) leads in the considered scheme. For example, in the distance-based scheme, the maximum can be calculated by $p_\mathrm{distance} = \max_{i\in I, j\in J} \{\nicefrac{C_{ij}^*}{d_i}\}$. Next, we identify optimal assignments for all compensation values $p\in \{ \nicefrac{\ell \cdot p_\mathrm{max}}{25}\mid \ell=0, 1, \dots 25\}$ by solving formulation~\eqref{eq:milp} using $	{\bar{c}}^*_{ij} = P_{ij}(C_{ij}) C_{ij} + (1 - P_{ij}(C_{ij})) c'_i$, where $C_{ij}$ is calculated using $p$ for all $i\in I$ and $j\in J$. Let $\ell^*\in \{0, 1, \dots, 25\}$ be a value yielding minimum expected costs for the resulting assignment. We use the golden section search within the interval $[\nicefrac{\ell^- \cdot p_{\mathrm{max}}}{25},\nicefrac{\ell^+ \cdot p_{\mathrm{max}}}{25}]$, $\ell^-= \max\{\ell^*-1,0\}$, $\ell^+=\min\{\ell^*+1,25\}$ to find a parameter value leading to a local minimum of the expected costs using formulation~\eqref{eq:milp}. This value is then used for the comparison of the different schemes in the following section.

\section{Computational study\label{sec:results}}

In this section, we discuss the results of our computational study for the \probabbr (with separability) using the exact two-stage solution approach proposed in \cref{sec:ilp}. We analyze the performance of the four considered compensation schemes for linear and logistic acceptance probability functions in \cref{sec:linearPerformance} and \cref{sec:logisticPerformance}, respectively. \cref{sec:sensitivity} analyzes the sensitivity of the results in the logistic model on the availability of occasional drivers, their willingness to take detours, and the penalty for rejecting tasks. \cref{sec:dynamic} analyzes the performance of our approach in a setting in which drivers and tasks appear dynamically over time.

We base our analysis of the compensation schemes on their respective performance in terms of economic benefit to the company
and satisfaction of the occasional drivers with the respective compensation schemes. We quantify these
performance indicators as follows:

\begin{itemize}
	\item To quantify the economic
	benefit, we consider the \textbf{expected total cost}, which is directly related to the objective function values of the \probabbr.
	\item To quantify the satisfaction of occasional drivers, we calculate the \textbf{mean acceptance rate} $\sum_{(i,j)\in A} \nicefrac{P_{ij}(C_{ij})}{|A|}$. A high value indicates that occasional drivers are likely to receive offers that they are satisfied with as compensation for completing the task. To put this indicator in perspective, we evaluate it together with the \textbf{fraction of tasks offered}. A high value for both indicators represents a compensation scheme in which 
	occasional drivers can expect to receive a high number of offers, and that these offers are generally acceptable.	
\end{itemize}

Our conclusions regarding comparisons between the different compensation schemes are supported by paired t-tests performed for each instance and reported at $\alpha=0.05$ level in the remainder of this section. Full results are provided in the electronic companion, which includes detailed results for each instance and performance indicator. Runtimes are not reported as they are less than one second for each instance. 

\subsection{Performance comparison for linear acceptance probabilities \label{sec:linearPerformance}}	

\begin{figure}[h]
\begin{subfigure}{0.32\textwidth}
	\includegraphics[width=\linewidth]{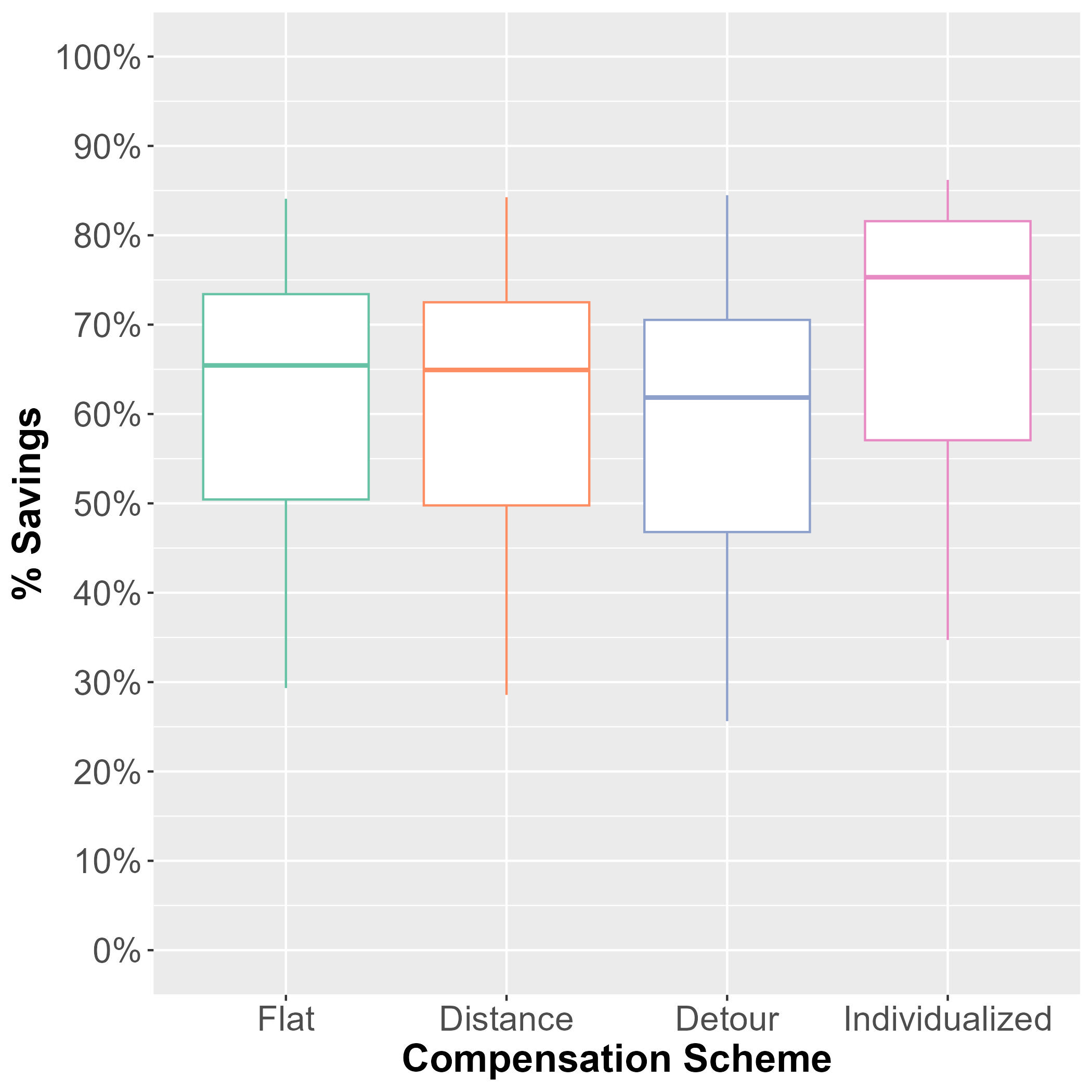}
	\caption{\label{fig:linearCost:boxplot}}
\end{subfigure}\hspace*{\fill}
\begin{subfigure}{0.64\textwidth}
	\includegraphics[width=\linewidth]{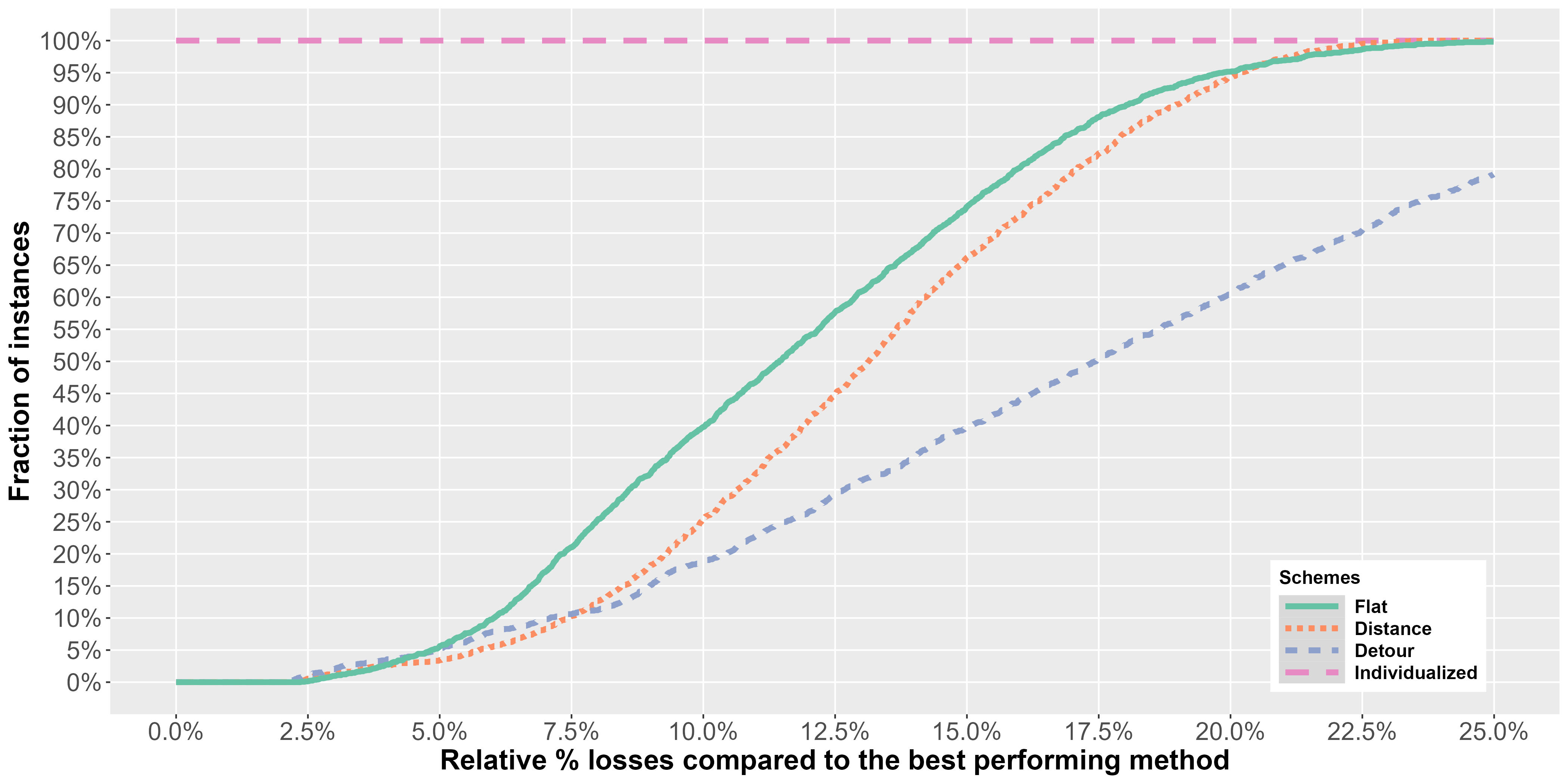}
	\caption{\label{fig:linearCost:performance}}
\end{subfigure}\hspace*{\fill}
\caption{Expected total cost for linear acceptance probabilities.} \label{fig:linearCost}
\end{figure}

\paragraph{Expected total cost} We focus first on analyzing the expected total cost. \cref{fig:linearCost:boxplot} shows the relative savings for this criterion for all compensation schemes compared to the case where all tasks are allocated to the 
professional drivers.
From this figure, we conclude that the individualized compensation scheme leads to cost savings with a median of over $75\%$. The other schemes lead to significantly smaller median savings of around $65\%$. This figure also indicates that the detour-based scheme performs the worst, while the results of the flat and distance-based schemes appear to be similar.

More insight into the relative performance of the four compensation schemes in terms of savings in expected total cost (in the linear case) can be obtained from \cref{fig:linearCost:performance}. For each scheme, this figure reports the fraction of instances in which the relative loss in solution quality (in percent) is at most a given value compared to the best performing scheme. 

First, we observe that the results obtained by using the individualized compensation scheme are strictly better than those of any other scheme in all instances considered. This follows because the relative loss of solution quality of the individualized scheme is zero in $100\%$ of the instances. For all other schemes, the relative loss is at least $2\%$ for each instance. \cref{fig:linearCost:performance} also illustrates that the individualized scheme clearly outperforms all other schemes, that the flat scheme performs slightly better than the distance scheme, and that both of them clearly outperform the detour scheme. Note that this relative order of the schemes is rather unexpected. Intuitively, one would expect the distance and detour schemes to outperform the flat scheme, since the detour and distance amounts directly affect the acceptance probabilities. 
%
The results show that for each experimental configuration, the intervals of the relative mean difference in expected total cost between the individualized and benchmarks schemes are $(4.32\%,79.51\%)$, $(6.10\%,61.28\%)$, and $(5.02\%,66.78\%)$ for the detour, distance, and flat schemes, respectively. Paired t-tests confirm that the differences between the compensation schemes are statistically significant in all instances.

\begin{figure}[h]	
\hspace*{\fill}
\begin{subfigure}{0.48\textwidth}
\includegraphics[width=\linewidth]{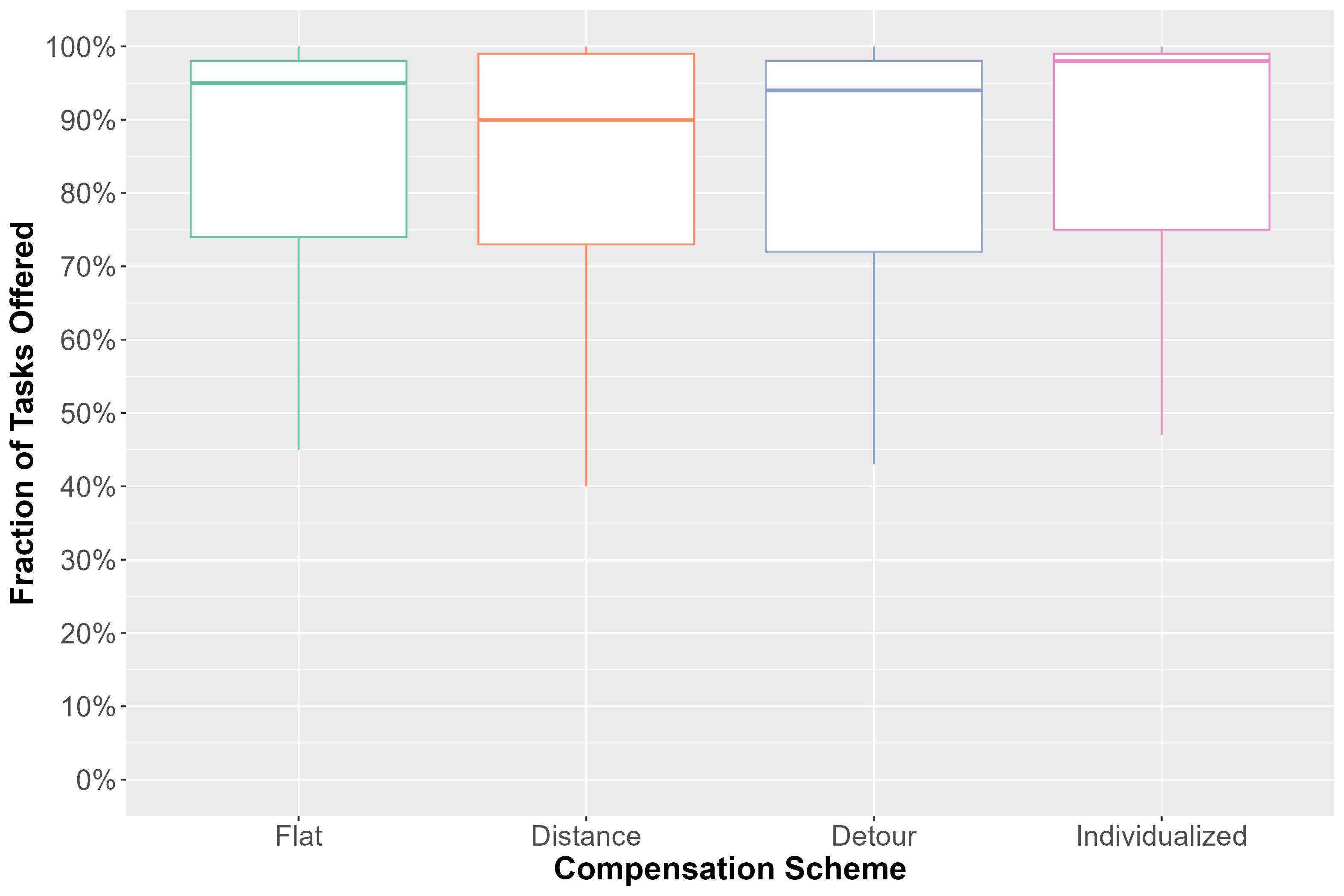}
\caption{} 
\label{fig:linearAcceptance:tasksoffered}
\end{subfigure}\hspace*{\fill}
\begin{subfigure}{0.48\textwidth}
\includegraphics[width=\linewidth]{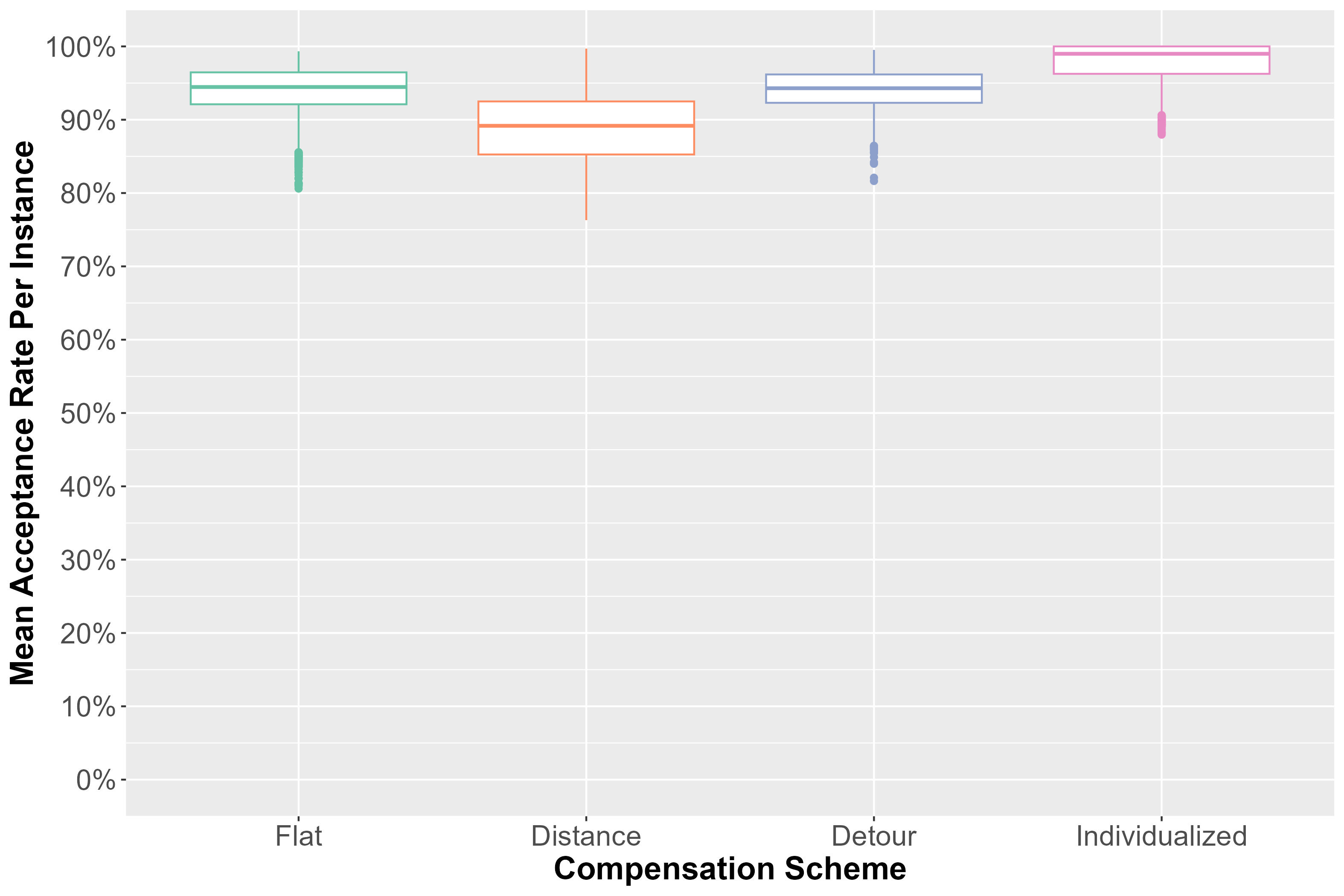}
\caption{}
\label{fig:linearAcceptance:meanacceptance}
\end{subfigure}\hspace*{\fill}
\caption{Fractions of offered tasks and expected acceptance rates for linear acceptance probabilities.}
\label{fig:linearAcceptance}
\end{figure}

\paragraph{Mean acceptance rate} 
\cref{fig:linearAcceptance:tasksoffered} shows the fractions of tasks offered to occasional drivers by the different compensation schemes,
while \cref{fig:linearAcceptance:meanacceptance} shows the mean acceptance rates of these tasks per instance. We observe that all schemes offer the majority of tasks to occasional drivers while achieving high expected acceptance rates. As for the previous two criteria, the individualized compensation scheme clearly outperforms the benchmark schemes as it achieves a higher acceptance rate while simultaneously offering more tasks. Among the benchmark schemes, the flat scheme shows the best performance, followed by the detour and distance schemes. 
The average differences in mean acceptance rates between the individualized and the benchmark schemes per experimental configuration take values in $(0.88\%,8.10\%)$, $(1.18\%,16.04\%)$ and $(1.24\%,8.10\%)$ for the detour, distance, and flat schemes, respectively. While the differences between the individualized scheme and the distance or flat scheme are statistically significant in all instances, this is not true for 49 out of 330 cases when comparing the individualized scheme and the detour scheme. These exceptions arise for instances with a relatively small number of occasional drivers, i.e., when $|J|\in \{50,75\}$. As their availability increases, the differences between the two schemes increase. 

Overall, we conclude that the individualized scheme clearly outperforms all alternatives considered in each
evaluation criteria. Thus, the use of the individualized scheme has the greatest potential to increase the satisfaction of occasional drivers and to reduce 
the use of professional drivers,
while simultaneously reducing the expected total cost. 

\subsection{Performance comparison for logistic acceptance probabilities\label{sec:logisticPerformance}}
\paragraph{Expected total cost} 

\begin{figure}[h]
\begin{subfigure}{0.32\textwidth}
\includegraphics[width=\linewidth]{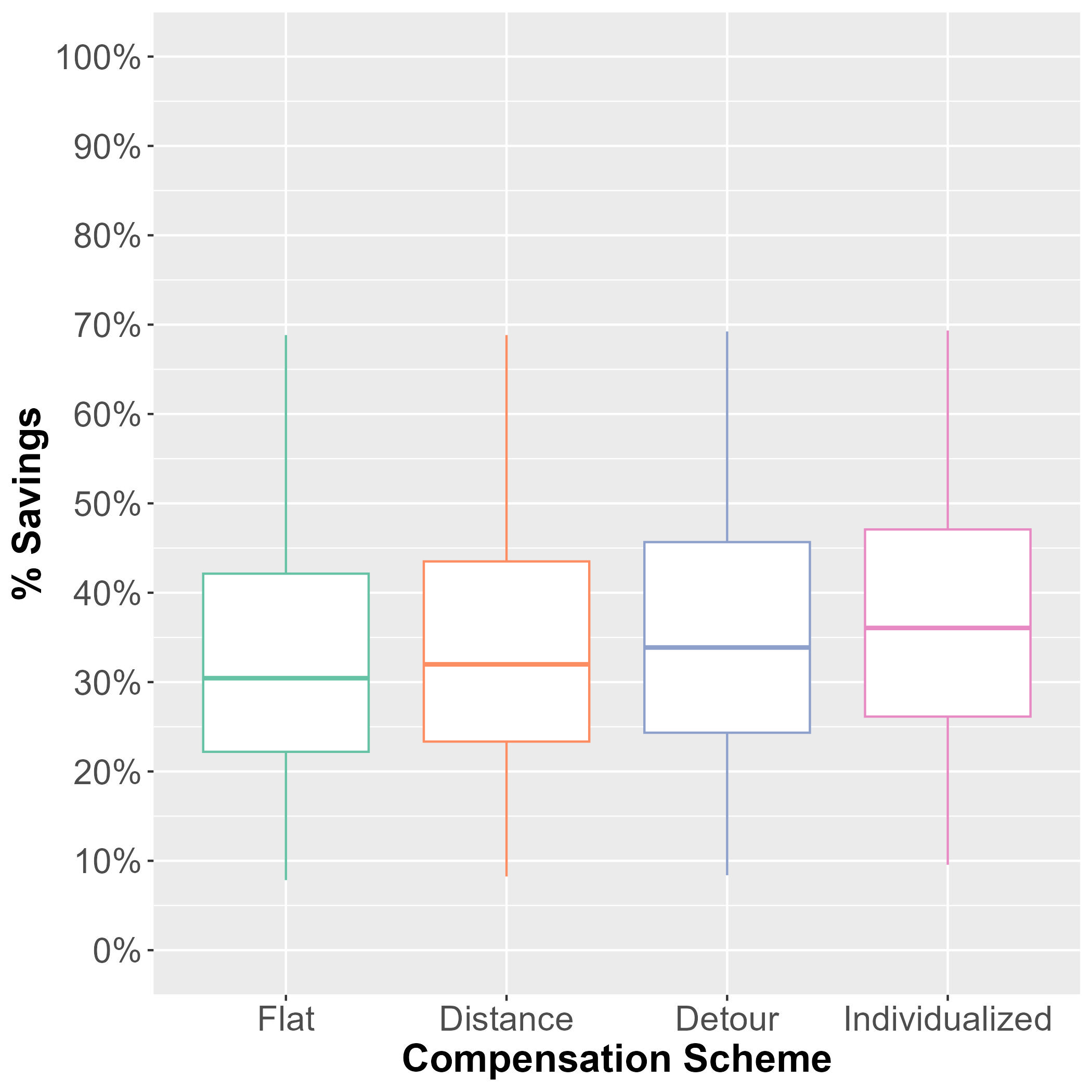}
\caption{}\label{fig:logisticCost:boxplot}
\end{subfigure}\hspace*{\fill}
\begin{subfigure}{0.64\textwidth}
\includegraphics[width=\linewidth]{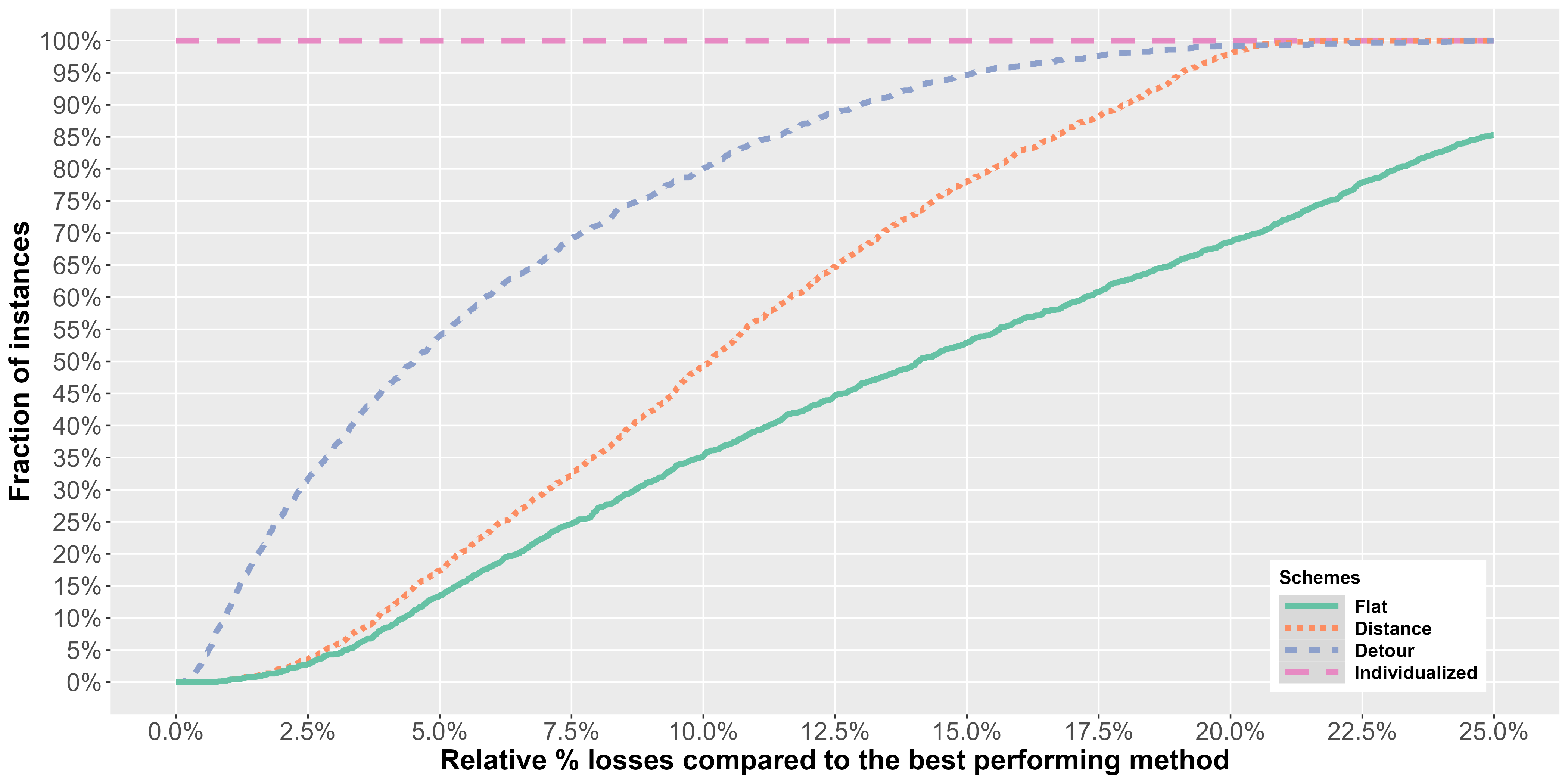}
\caption{}\label{fig:logisticCost:performance}
\end{subfigure}
\caption{Expected total cost for logistic acceptance probabilities.}
\label{fig:logisticCost}
\end{figure}

\cref{fig:logisticCost:boxplot} shows the relative cost savings of all four compensation schemes compared to the setting with no occasional drivers for logistic acceptance probability functions. Consistent with the case of the linear acceptance function, the individualized compensation scheme outperforms all benchmark schemes with median cost savings of more than $35\%$. \cref{fig:logisticCost:performance} shows that the individualized scheme outperforms all other schemes in every instance. For each parameter combination, the mean differences between the individualized scheme and the detour-based, distance-based, and flat schemes take values in the intervals $(0.27\%,6.73\%)$, $(1.31\%,8.99\%)$ and $(1.41\%,13.60\%)$, respectively. The statistical significance of the differences between the individualized and the benchmark schemes is confirmed for each setting considered by paired t-tests. 
\cref{fig:logisticCost} also reveals that the detour-based scheme comes closest to the individualized scheme and outperforms the other two schemes. The flat compensation scheme shows the worst performance, which is in contrast to the case of linear acceptance probability functions, where the performance order of the three benchmark schemes is reversed.

\begin{figure}[h]
\hspace*{\fill}
\begin{subfigure}{0.48\textwidth}
\includegraphics[width=\linewidth]{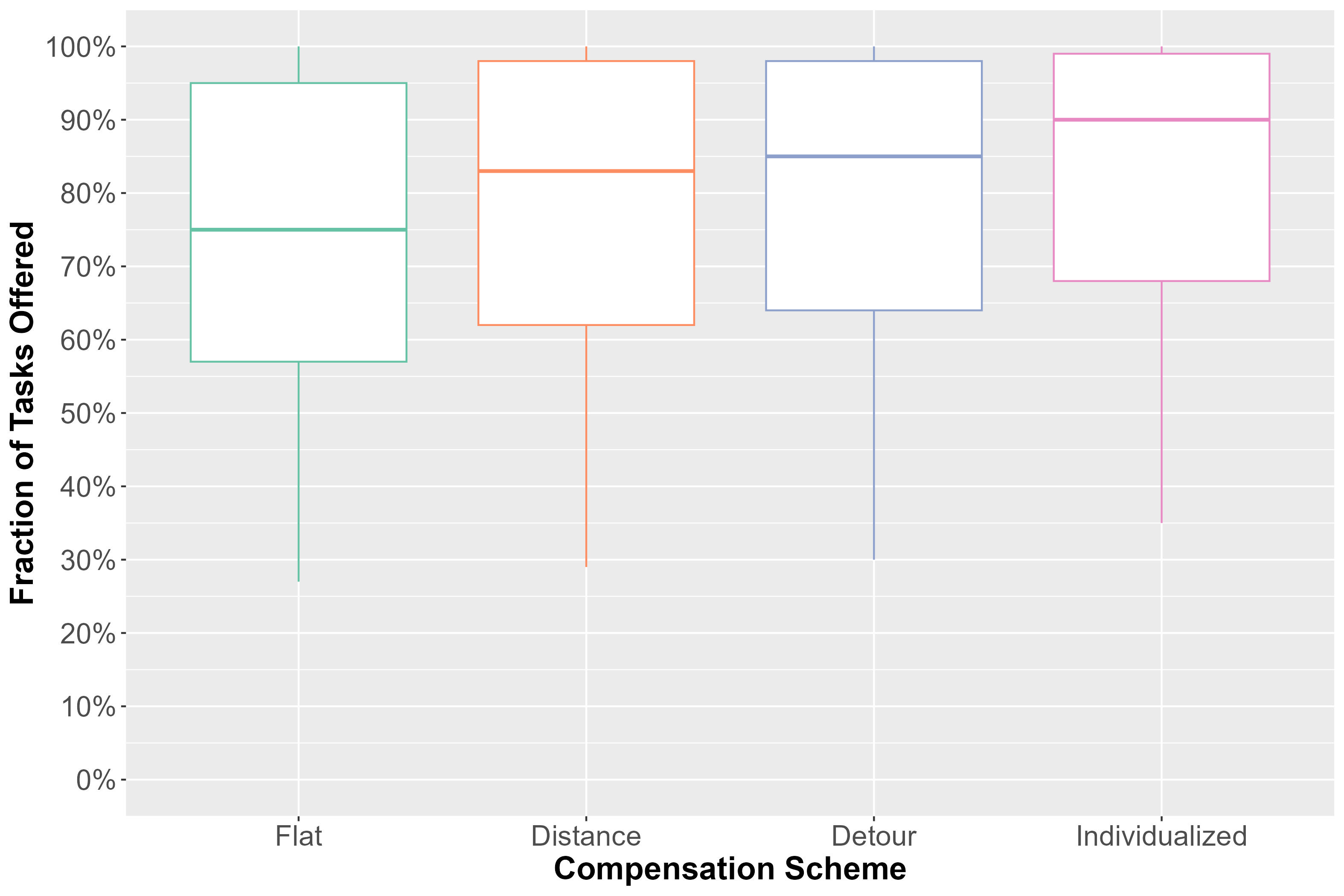}
\caption{}\label{fig:logisticAcceptance:tasksoffered}
\end{subfigure}\hspace*{\fill}
\begin{subfigure}{0.48\textwidth}
\includegraphics[width=\linewidth]{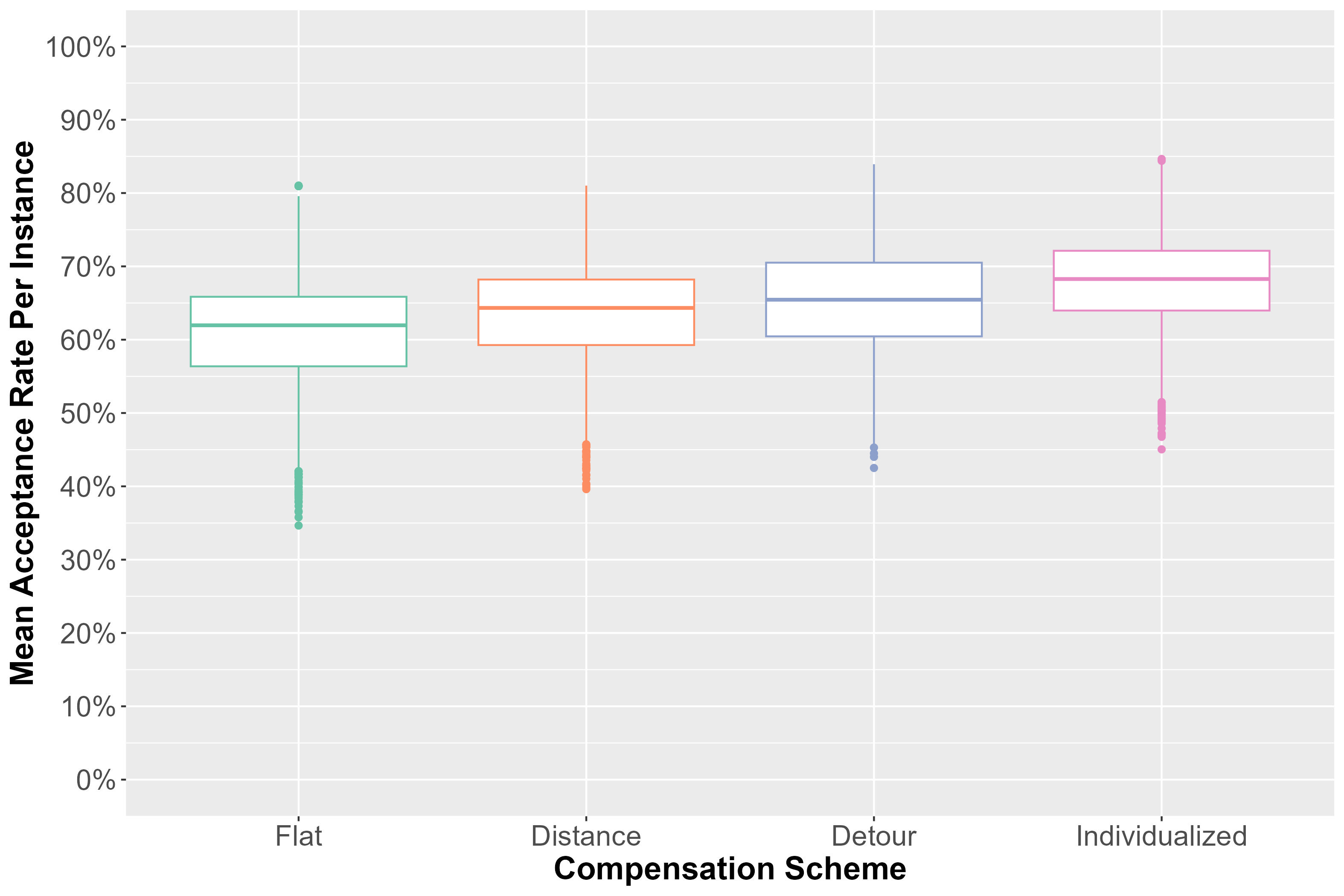}
\caption{}\label{fig:logisticAcceptance:meanacceptance}
\end{subfigure}\hspace*{\fill}
\caption{Fractions of offered tasks and expected acceptance rates for logistic acceptance probabilities.}\label{fig:logisticAcceptance}
\end{figure}

\paragraph{Mean acceptance rate}
\cref{fig:logisticAcceptance:tasksoffered,fig:logisticAcceptance:meanacceptance} show that the individualized scheme is superior to the benchmark schemes, as it outsources more tasks than the other schemes and achieves the highest mean acceptance rate. The median fraction of tasks offered to occasional drivers is around $90\%$, and the median expected acceptance rate (per instance) of these offers is close to $70\%$.
In comparison, the median fraction of tasks offered by the flat compensation scheme is only about $75\%$, and the median acceptance rate per instance is only slightly higher than $60\%$. If we examine the differences in each configuration, we see that the average differences of the mean acceptance rates per configuration of the detour, distance, and flat schemes to the individualized scheme lie in the intervals $(1.04\%,11.75\%)$, $(1.30\%,14.13\%)$, and $(3.11\%,24.68\%)$, respectively. While the differences between the individualized and the flat schemes are statistically significant in all cases, this is not true for the distance scheme in 2 out of $330$ configurations ($|J|=150$, $\rho\in \{20,50\}$, $\mu=0$) in which occasional drivers are mainly concerned with distance and are not sensitive to the detour. When comparing the detour-based and the individualized scheme, we observe that the differences are not statistically significant in $76$ out of $330$ configurations (in which $|J|\in \{50,75,100\}$). This shows that the individualized scheme significantly outperforms the detour-based schemes in environments characterized by an oversupply of occasional drivers.

Overall, we conclude that the individualized scheme clearly outperforms all alternatives considered in each
evaluation criteria for both linear and logistic acceptance probability functions. 

\subsection{Sensitivity analysis} \label{sec:sensitivity}
In the following, we analyze the sensitivity of the
performance indicators with regard to changes in the availability of occasional drivers $|J|$, the penalty for rejected tasks $\pnlty$, and the 
parameter $\utlty$ that affects the willingness to make detours. We also analyze the impact of inaccuracies in the prediction of acceptance probabilities. For the sake of brevity, and since all effects can be readily demonstrated within the logistic acceptance probability model, we restrict our discussion of sensitivity to the above parameters to that model. The full set of results can be found in the electronic companion. 

\begin{figure}[h]
\hspace*{\fill}
\begin{subfigure}{0.32\textwidth}
\includegraphics[width=\linewidth]{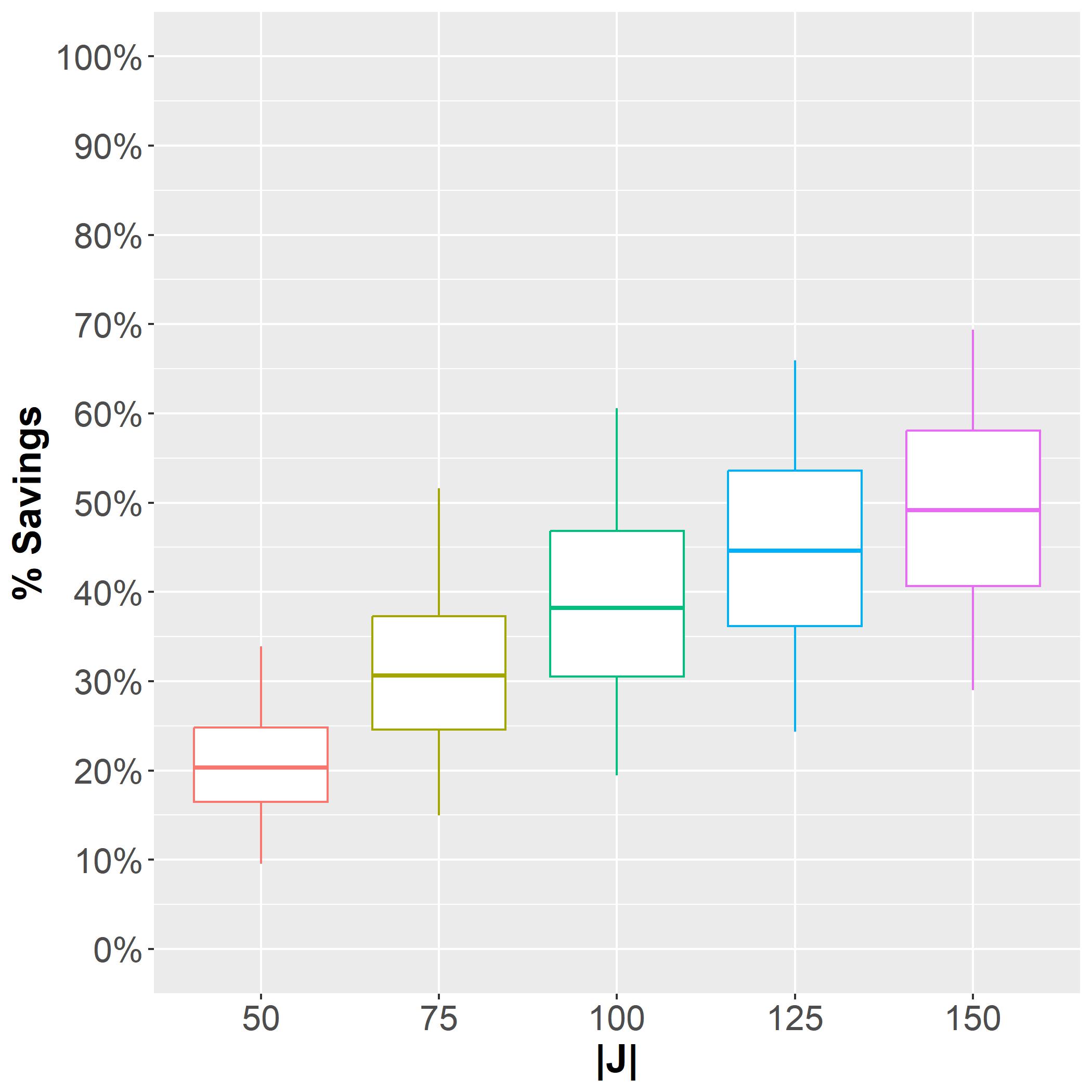}
\caption{Expected total cost.}\label{fig:logistic_cost_O}
\end{subfigure}\hspace*{\fill}
\begin{subfigure}{0.32\textwidth}
\includegraphics[width=\linewidth]{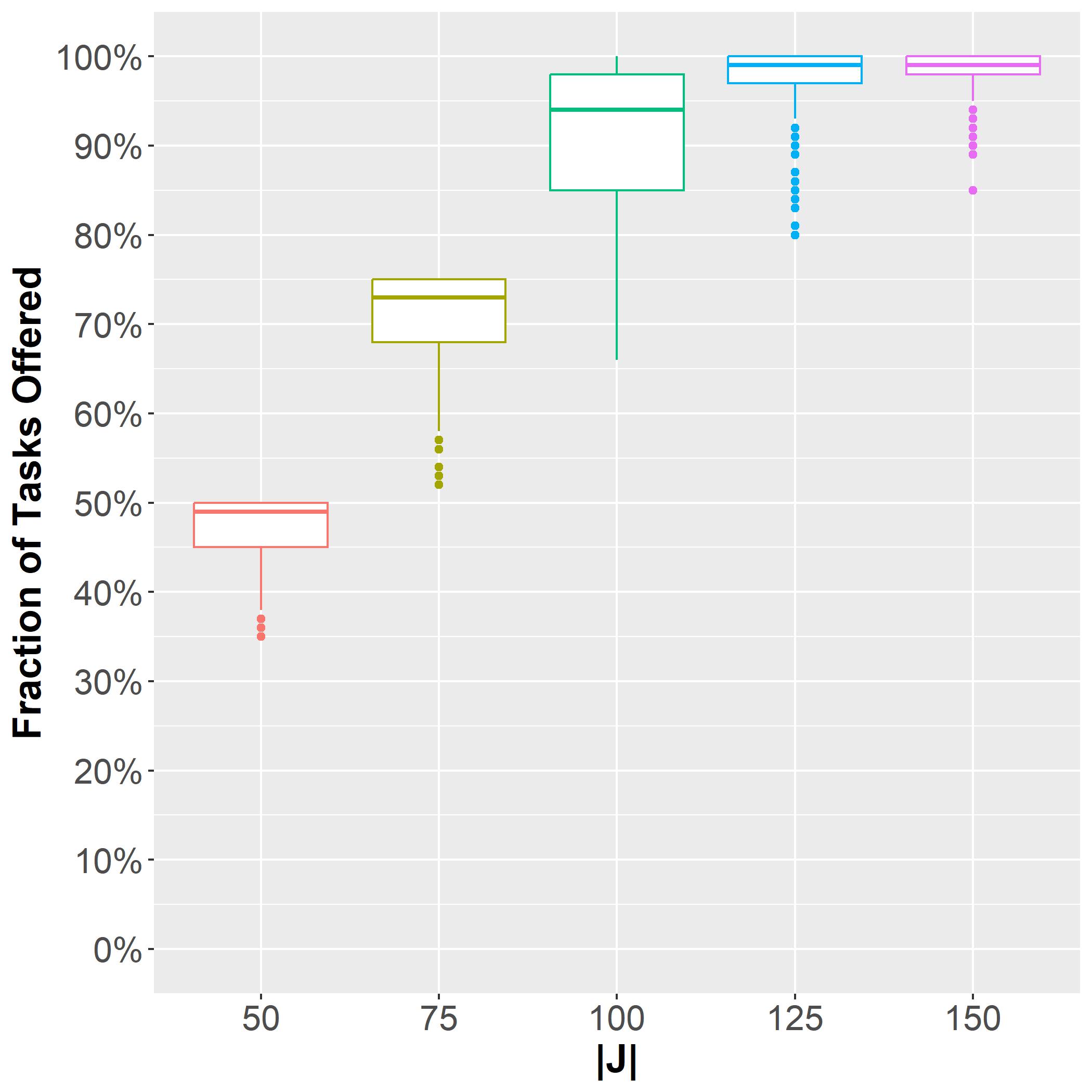}
\caption{Fraction of tasks offered.}\label{fig:logistic_frac_acc_O}
\end{subfigure}\hspace*{\fill}
\begin{subfigure}{0.32\textwidth}
\includegraphics[width=\linewidth]{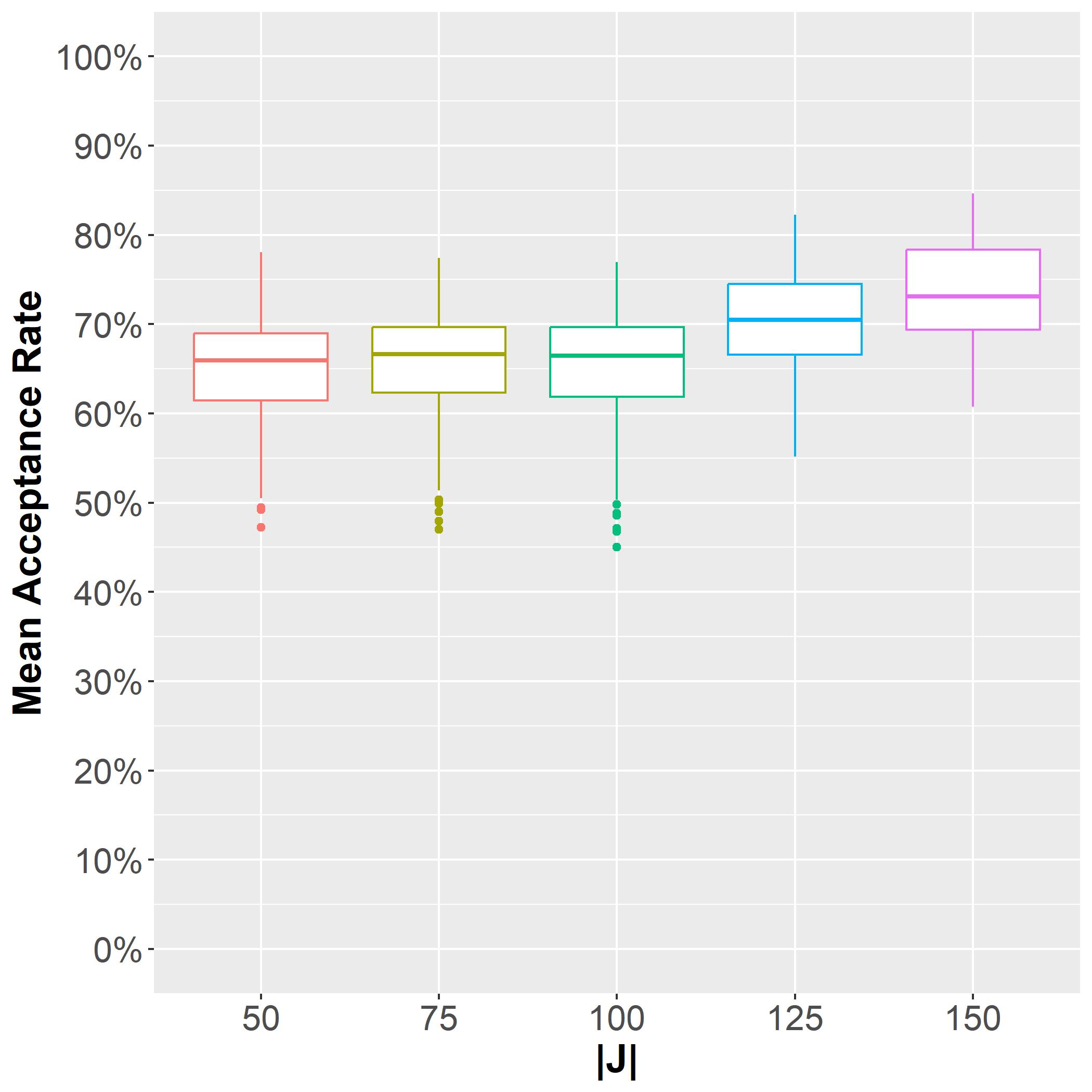}
\caption{Mean acceptance rate.}\label{fig:logistic_mean_acc_O}
\end{subfigure}\hspace*{\fill}		
\caption{Sensitivity to availability} \label{fig:sens_od}
\end{figure}

\paragraph{Sensitivity towards availability of occasional drivers} \cref{fig:sens_od} shows the savings in expected total cost
compared to the scenario with no occasional drivers, as well as the mean acceptance rates and fraction of tasks offered, respectively, for different numbers of (available) occasional drivers $|J|$. We observe that the savings tend to increase with increasing availability of occasional drivers for all the performance indicators, which is to be expected since a higher number gives the company more options to select cost-efficient occasional drivers. 
\cref{fig:logistic_frac_acc_O} nicely illustrates that the number of tasks offered is always close to the maximum number possible (i.e., the minimum of $|I|$ and $|J$), while the mean acceptance always hovers around $65\%$. With a higher number of occasional drivers (i.e., with $|J| > |I|$), it is more likely to find occasional drivers who are more willing to accept the offers, and the mean acceptance increases (cf.\ \cref{fig:logistic_mean_acc_O}).

\begin{figure}[h]
\hspace*{\fill}
\begin{subfigure}{0.32\textwidth}
\includegraphics[width=\linewidth]{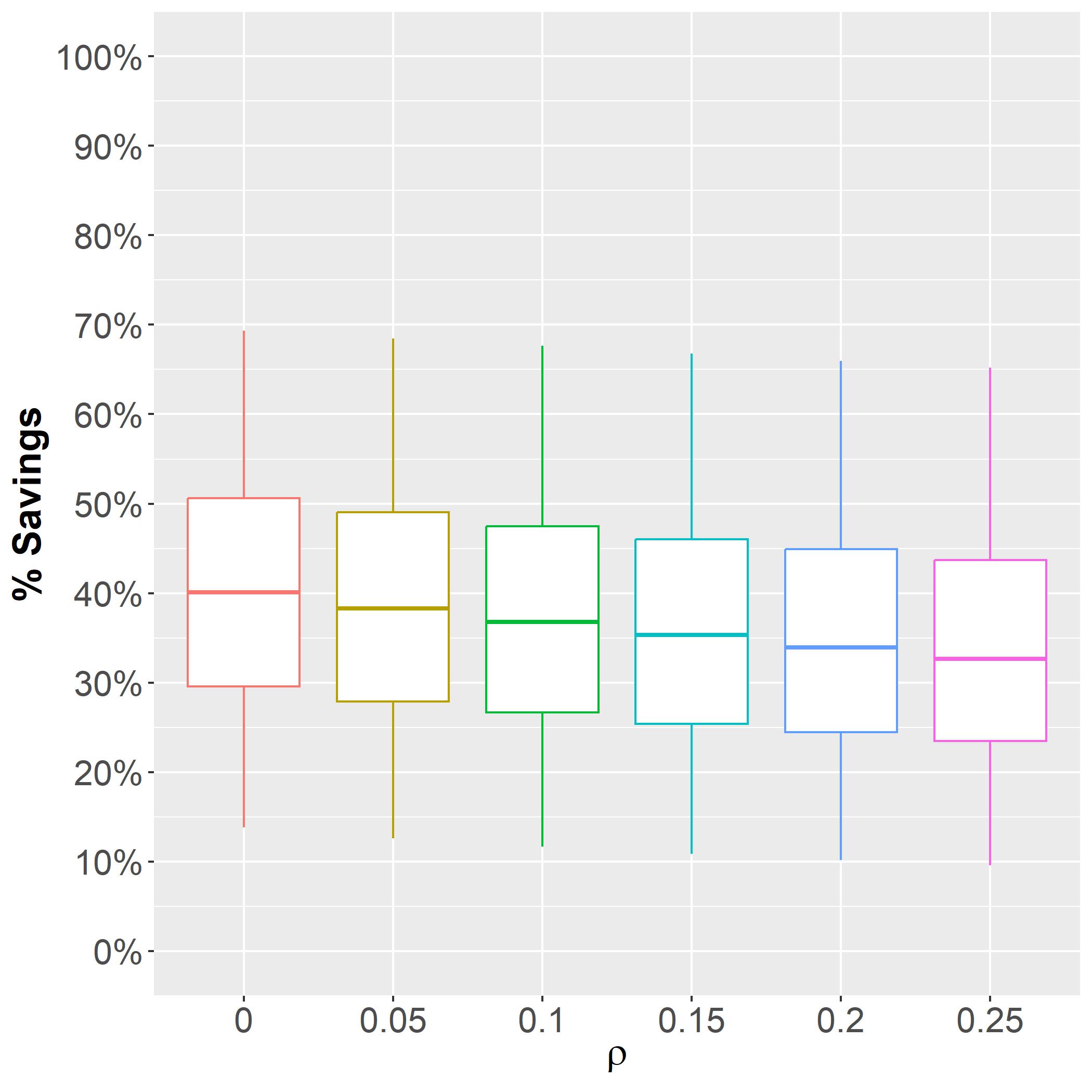}
\caption{Expected total cost.}\label{fig:logistic_cost_P}
\end{subfigure}\hspace*{\fill}
\begin{subfigure}{0.32\textwidth}
\includegraphics[width=\linewidth]{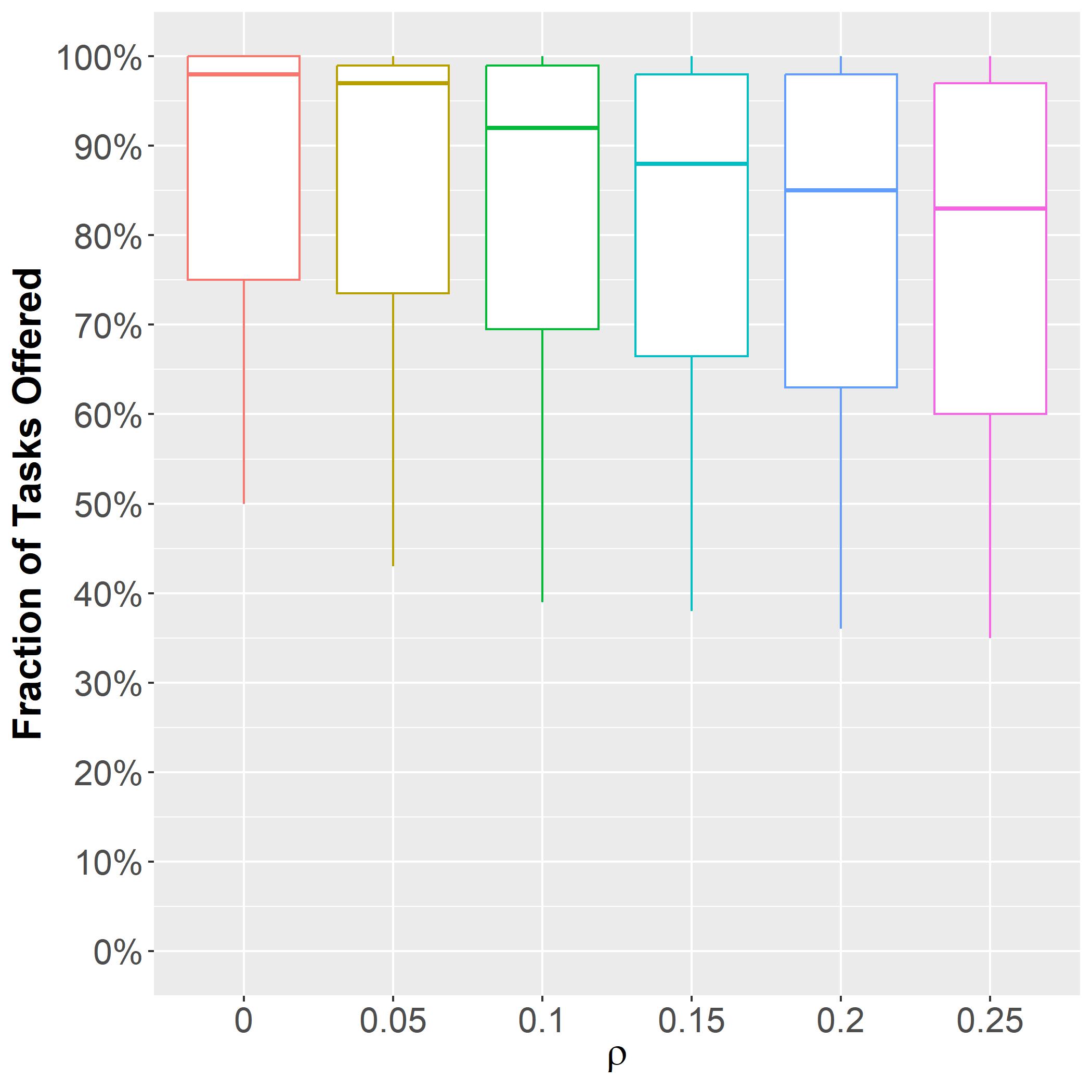}
\caption{Fraction of tasks offered.}\label{fig:logistic_frac_acc_P}	
\end{subfigure}\hspace*{\fill}			
\begin{subfigure}{0.32\textwidth}
\includegraphics[width=\linewidth]{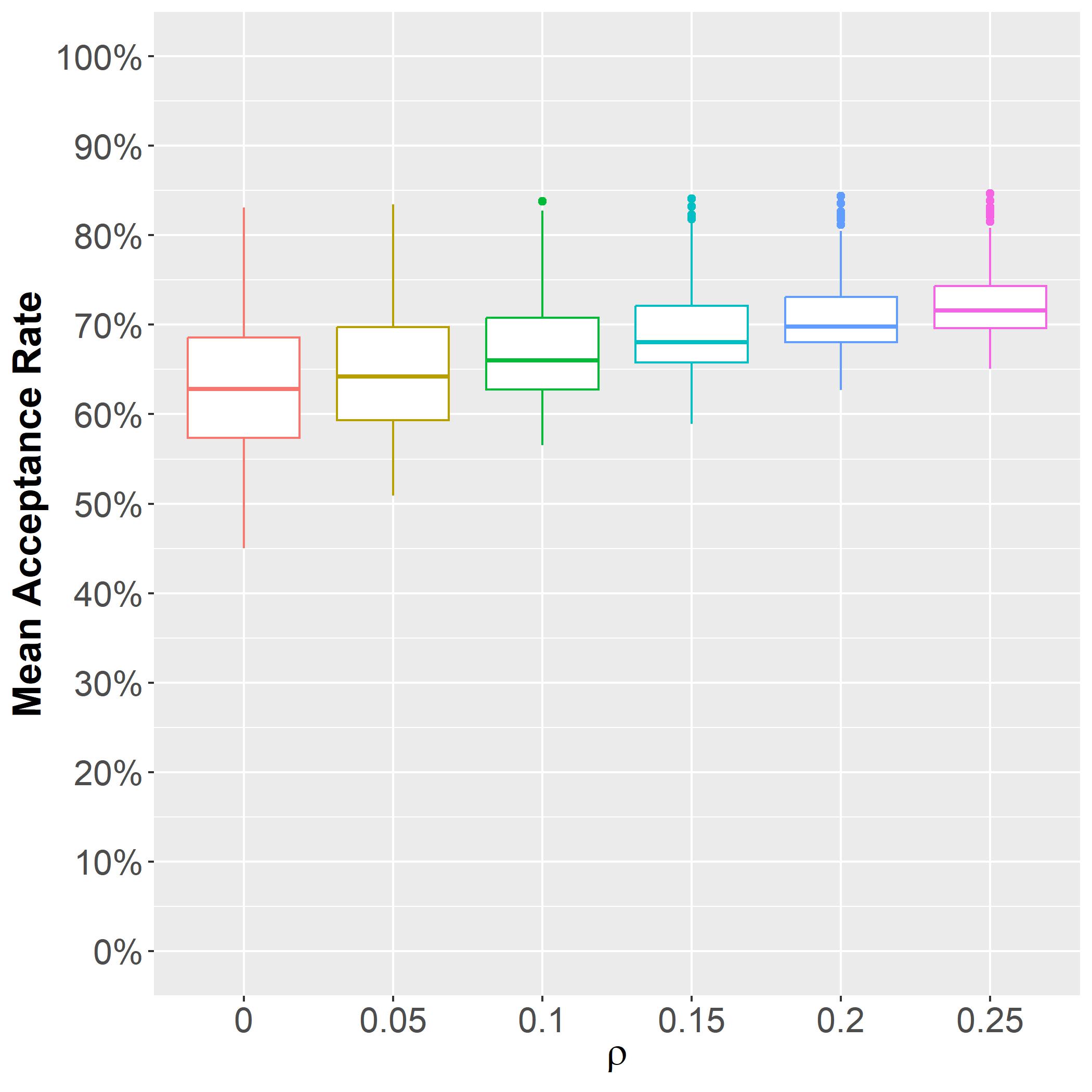}
\caption{Mean acceptance rate.}\label{fig:logistic_mean_acc_P}
\end{subfigure}\hspace*{\fill}			
\caption{Sensitivity to penalty}\label{fig:sens_pen}
\end{figure}

\paragraph{Sensitivity towards the penalty parameter $\pnlty$} \cref{fig:sens_pen} illustrates the effect of increasing the penalty for rejected offers. As expected, an increase in the penalty has a significant impact on the individualized scheme, leading to fewer tasks being offered to the occasional drivers on average, cf.\ \cref{fig:logistic_frac_acc_P}. At the same time, the mean acceptance of offers increases (cf.\ \cref{fig:logistic_mean_acc_P}), indicating a more careful selection of occasional drivers and a stronger focus on making good (acceptable) offers. \cref{fig:logistic_cost_P} shows that the decrease in the number of occasional drivers used and the simultaneous increase in the penalty also lead to smaller savings in terms of total expected costs. However, as can be seen in \cref{fig:logistic_cost_P} as well, the effects are much less severe than one would expect, decreasing from $40\%$ savings in the unrealistic case that no penalty is incurred to $33\%$ in the case that a rejected offer leads to a $25\%$ increase in operational costs. One explanation is that the individualized compensation scheme allows good (acceptable) offers to be made to occasional drivers for tasks that would initially be very costly to the company using professional drivers, thus allowing the company to efficiently outsource them to occasional drivers. 

\begin{figure}[h]	
\hspace*{\fill}
\begin{subfigure}{0.32\textwidth}
\includegraphics[width=\linewidth]{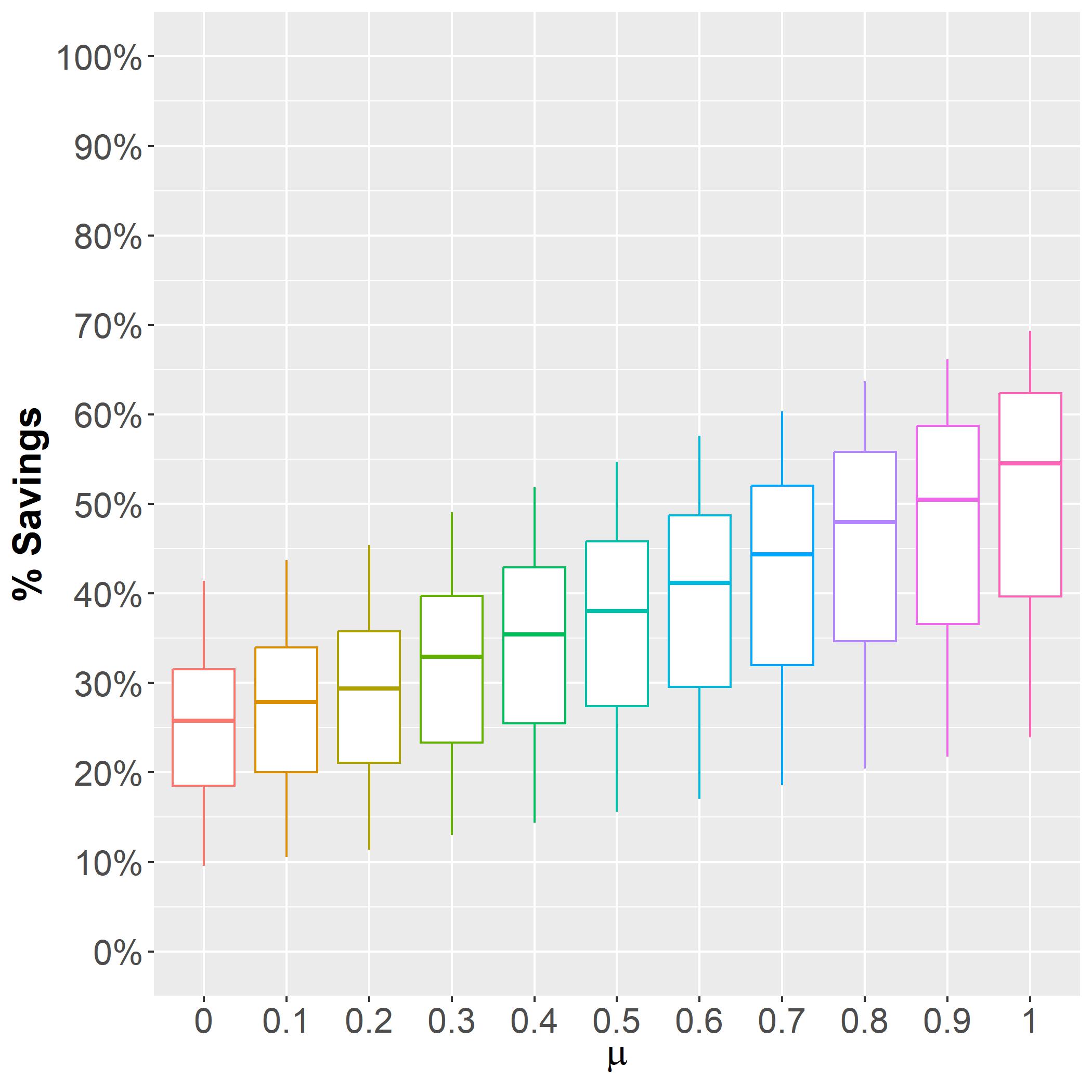}
\caption{Expected total cost.}\label{fig:logistic_cost_R}
\end{subfigure}\hspace*{\fill}
\begin{subfigure}{0.32\textwidth}
\includegraphics[width=\linewidth]{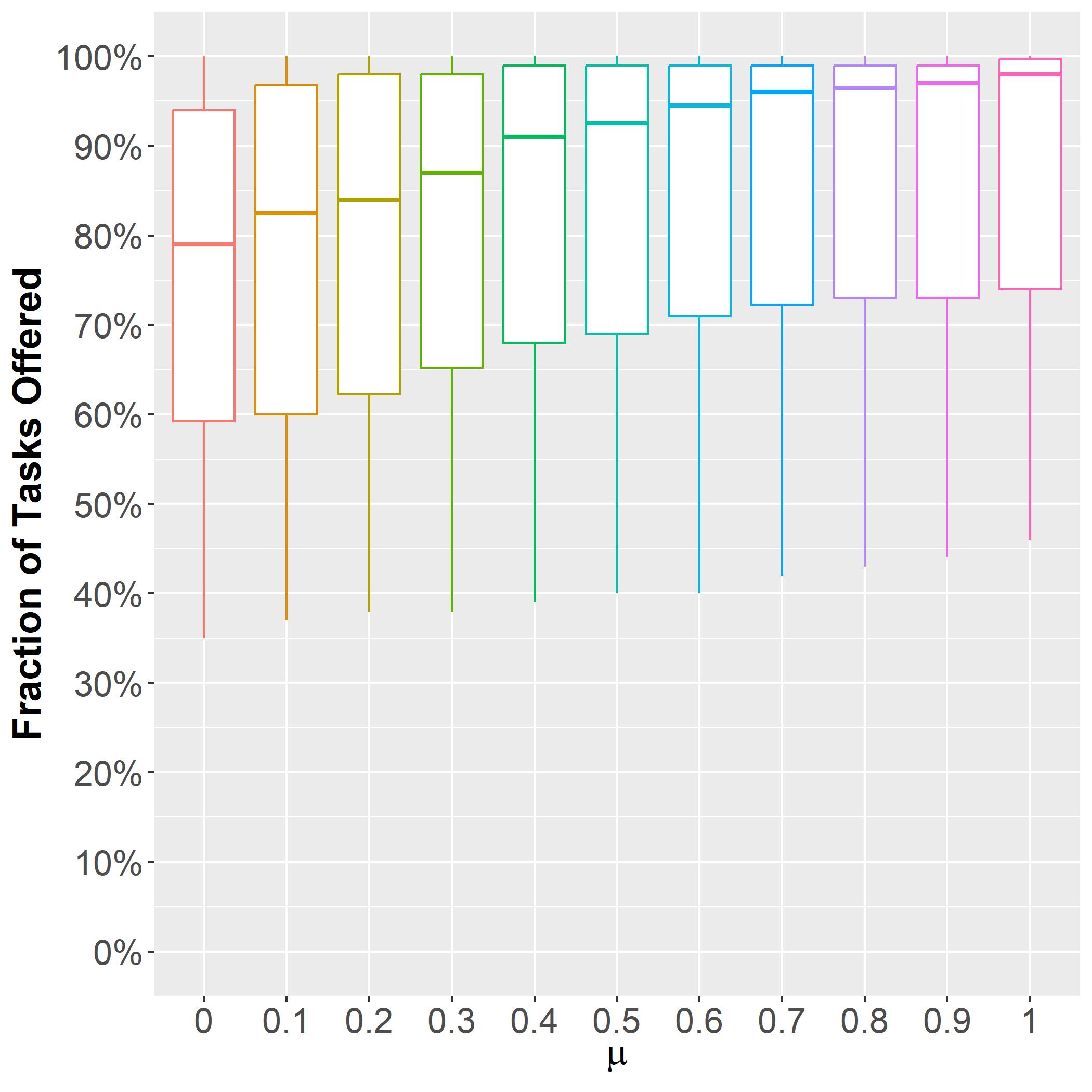}
\caption{Fraction of tasks offered.}\label{fig:logistic_frac_acc_R}
\end{subfigure}\hspace*{\fill}
\begin{subfigure}{0.32\textwidth}
\includegraphics[width=\linewidth]{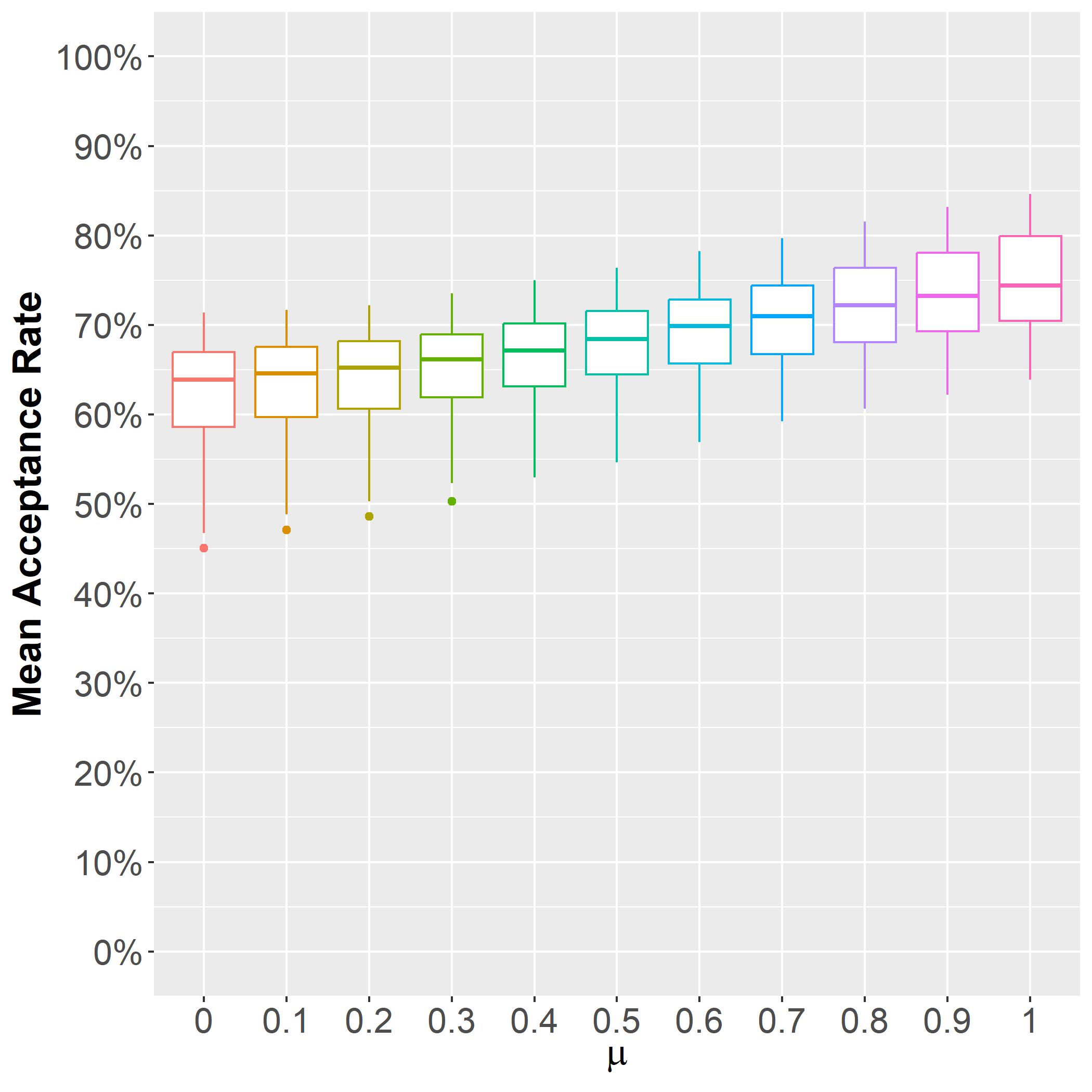}
\caption{Mean acceptance rate.}\label{fig:logistic_mean_acc_R}
\end{subfigure}\hspace*{\fill}

\caption{Sensitivity to distance utility weight} \label{fig:sens_util}
\end{figure}

\paragraph{Sensitivity towards parameter $\utlty$} A larger value of $\mu$ increases the acceptance probability because it increases the willingness of occasional drivers to make a detour. As can be seen in \cref{fig:sens_util}, this increases 
the savings in terms of expected total cost,
more tasks are offered to occasional drivers, and the acceptance rates also increase, cf.\ \cref{fig:logistic_cost_R,fig:logistic_frac_acc_R,fig:logistic_mean_acc_R}.
While this effect was expected, it is worth noting that the individualized scheme generates significant cost
savings for all parameter values considered. Furthermore, high fractions of tasks offered to occasional drivers and mean acceptance rates indicate that it allows to outsource a high proportion of tasks to occasional drivers, relatively independent of the concrete value of parameter $\utlty$. 

\begin{figure}[h]
\centering
\includegraphics[width=0.6\textwidth]{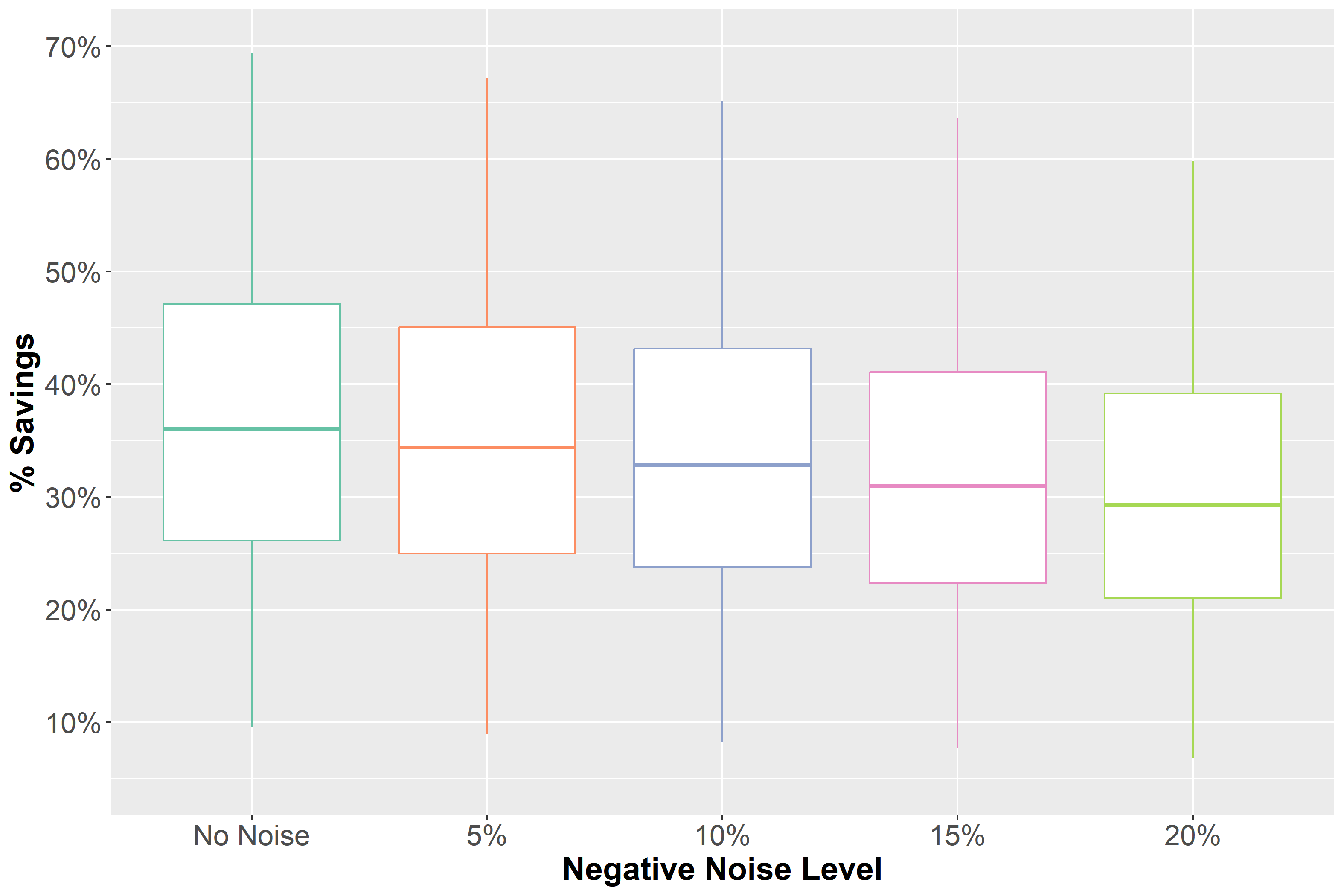}
\caption{Sensitivity to model accuracy} \label{fig:accuracy}	
\end{figure}

\paragraph{Sensitivity towards acceptance probability model inaccuracy} 
To demonstrate the impact of prediction inaccuracy on the expected cost savings, we introduce an additive noise term 
that represents the difference 
between the real acceptance probabilities and the estimated ones used by the model to decide on offers to crowdshippers and associated compensation values. 
In our analysis, we focus on scenarios with negative noise that
decreases the acceptance probability and leads to larger expected costs. 
For every task-driver pair in the solution of an instance and each considered noise level $n\in \{0.05, 0.1, 0.15, 0.2\}$, we generate a negative noise uniformly at random from the interval $[-n,0]$. For each noise level and every solution, we then use the noised predictions to recalculate the expected cost savings compared to the base scenario in which no occasional drivers are utilized. 
\cref{fig:accuracy} shows the distribution of expected cost savings. We observe that significant expected cost savings can be achieved even with increasing noise levels. 
Furthermore, increasing noise level leads to a reduction in expected cost savings, which drops from a median of $36\%$ to $29\%$ when the noise parameter is increased from $0$ to $0.2$. We conclude that the more accurate the acceptance probability models used within the proposed solution method, the higher the potential of obtaining solutions with larger expected cost savings.

\subsection{Performance in a dynamic environment} \label{sec:dynamic}

In many cases, tasks and drivers appear dynamically over time and sequential decisions 
about the allocation of tasks and compensations offered are required
\citep{allahviranloo2019dynamic,ausseil2022supplier}. Unfortunately, the \probabbr loses some of its strong theoretical properties in this setting and we can no longer find the provably optimal compensation according to the acceptance probability of drivers. 
The main reason for this is that 
defining values for the costs $c$ and $c'$ becomes difficult
in a dynamic setting since not assigning a task within one decision period might not rule out assigning the task in a later period.
Despite these difficulties, we will show in the following, that the \probabbr can be employed as a powerful heuristic in the dynamic setting as well. In the following, we first describe the dynamic environment in which we then 
benchmark the \probabbr against other compensation schemes from the literature.

\begin{figure}[h]
\centering
\resizebox{0.7\textwidth}{!}{%
\begin{tikzpicture}[node distance = 2cm, auto]
\node [block] (start) {Start};
\node [decision, below of=start, xshift=3cm] (compare) {Current period $\leq$ Final period};
\node [block, below of=start, node distance = 6cm] (end) {Finish Simulation};
\node [block, right of=start, node distance = 7cm] (remove) {Remove tasks and drivers that reached deadline};
\node [block, right of=remove, node distance = 3cm ] (update) {Generate new drivers and tasks};
\node [block, right of=update, node distance = 3cm ] (generate) {Generate offers};
\node [block, below of=generate, node distance = 6cm] (simulate) {Simulate driver decisions};
\node [block, below of=update, node distance = 6cm] (assign) {Remove drivers and tasks for all accepted offers};	
\node [block, below of=remove, node distance = 6cm] (updateperiod) {Update current period};

\path [line] (start)  |- (compare);
\path [line] (compare) |- node {Yes} (remove);
\path [line] (remove) -- (update);
\path [line] (update) -- (generate);
\path [line] (generate) -- (simulate);
\path [line] (simulate) -- (assign);
\path [line] (assign) -- (updateperiod);
\path [line] (updateperiod) -| ++(-2,0) |- (compare);
\path [line] (compare) |- node {No} (end);
\end{tikzpicture}
}
\caption{Flowchart of simulation mechanism}\label{fig:dynamic:flowchart}
\end{figure}
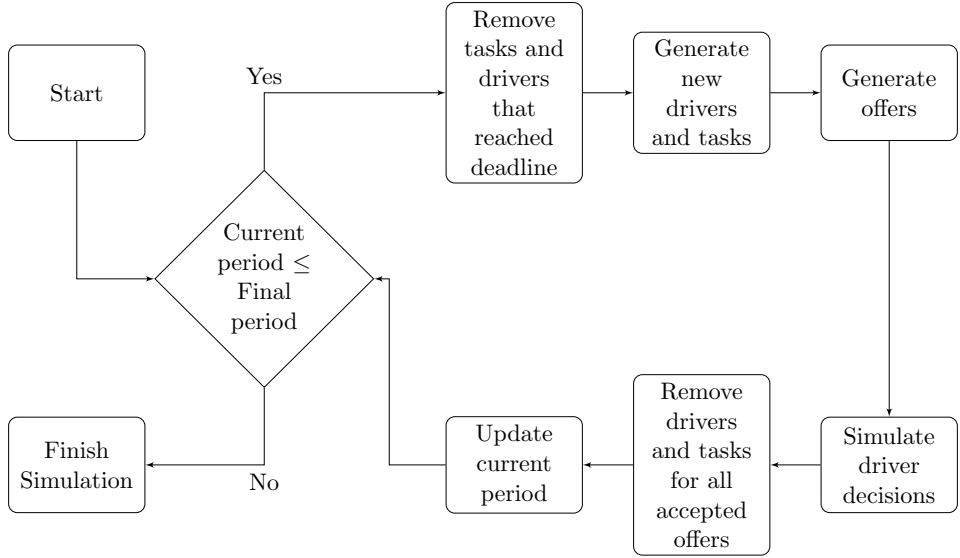

To evaluate the \probabbr in the dynamic setting, we simulate an environment in which the available drivers and tasks are revealed progressively over time. The process is described in Figure \ref{fig:dynamic:flowchart}. The planning horizon is discretized into periods, and new tasks and drivers may appear at the beginning of each period. Drivers and tasks stay in the system for a certain number of periods, and we assume that these numbers are known when the driver or the task is revealed. 
At each period the planner decides which offers (consisting of a task and offered compensation) to make to available drivers. As in the static version of the problem, drivers can accept or reject offers. 
If a driver accepts an offer, the driver and the task leave the system. Otherwise, they stay in the system for the next period unless they have reached their end period. 
An unassigned task leaving the system can either be considered as a lost sale or as outsourced to a 3PL as in the static setting. We assume the planner does not have any information about the availability of future drivers and tasks.

We generate instances for the dynamic case 
by varying the task and driver arrival rates as well as 
the maximum times that tasks and drivers stay in the system. Task and driver arrival rates are generated using a Poisson distribution with means $\mu \in \left\{5,10,15\right\}$. The maximum number of periods they remain in the system are generated using a Poisson distribution with means $\mu \in \left\{2,3,4\right\}$. For each combination of the means, we generate 25 realizations overall leading to 2\,025 problem instances. 
We solve each instance 25 times to account for the uncertainty in driver behavior (i.e., probabilistic acceptance and rejection of offers) leading to a total of 50\,625 computational experiments for each compensation scheme.
The planning horizon is partitioned in 100 periods.


For determining the compensation offered in the individualized scheme, we solve an instance of the \probabbr at each period with the current drivers and tasks. One difficulty lies in defining the parameters $c$ and $c'$ for the current set of tasks.
Indeed, these costs do only materialize in the last period a task is active and only if tasks not covered by crowdshippers are finally outsourced to a 3PL. 
In our experiments for the dynamic case we do, however, consider tasks not performed by crowdshippers as lost sales.
Therefore, we set both costs (i.e., $c$ and $c'$) to the distance of the task from the depot when the task appears and increase them linearly up to two times the distance in the final period. This ensures that the individualized scheme prioritizes tasks that are close to their end period and always tries to assign tasks. For the distance scheme, we offer a compensation equal to half of the distance of the task from the depot while the detour offers compensation equal to the detour the driver has to make to complete the task. 
It is more difficult to find a good compensation value for the flat scheme which does not consider task or driver related properties. Therefore, for the flat scheme we report results for three different compensation values.

\begin{figure}[h]	
\centering
\includegraphics[width=0.6\linewidth]{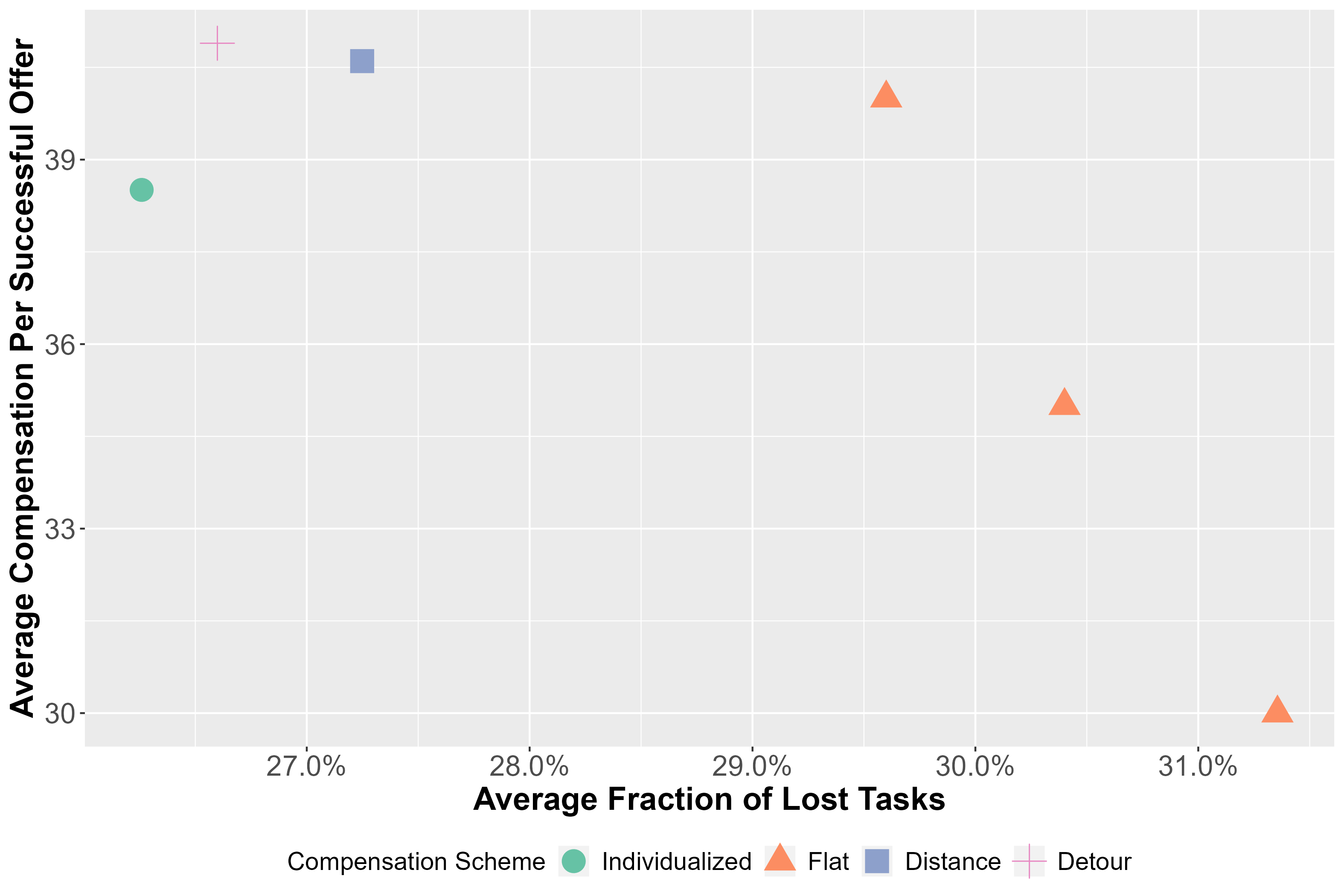}
\caption{Average performance of the schemes in terms of the average compensation per successful offer and the average fraction of lost tasks.}\label{fig:dynamic::simple}
\end{figure}

The results are summarized in Figure \ref{fig:dynamic::simple}. Each scheme is positioned based on the average performance in terms of the total number of assigned tasks and the compensation cost. The vertical axis shows the average performance of the schemes in terms of average compensation per successful offer 
and the horizontal axis the average fraction of lost tasks.

We see that the individualized scheme dominates the detour and distance schemes in terms of average compensation while achieving better performance in terms of lost tasks. The flat scheme is the worst-performing one in terms of the number of lost tasks but the best-performing in terms of average compensation for all three compensation values considered. 

The results indicate that the \probabbr can bring benefits even when used in a dynamic environment. The performance of the individualized scheme is, however, not as dominant as in the static setting. The reason could be that, in a dynamic environment, there are multiple factors under uncertainty not taken into account in any of the schemes (e.g., drivers and tasks appearing time) which may level out the performance of the schemes. Although this topic deserves more investigation, it goes beyond the scope of this analysis and it is left for future research.

\section{Conclusions and future work} \label{sec:conclusion}
In this paper, we introduced the \probname which defines an integrated solution to the task assignment and compensation decisions, while explicitly accounting for the probability with which occasional drivers will accept these assignments. To this end, we proposed an MINLP formulation for generic acceptance probability functions and showed that it can be solved optimally with a two-stage approach that decomposes compensation and assignment decisions. Moreover, we proposed an exact linearization of our MINLP formulation and showed that the latter can be solved in polynomial time for generic probability functions under mild conditions. We also derived explicit formulas for optimal compensation decisions for linear and logistic acceptance probability functions which are two particularly relevant, special cases of our generic setting.
Furthermore, we investigated the non-separable case where the decision-maker has additional restrictions on compensation and assignment decisions, and propose an approximate linearization for this extension.

We conducted an extensive computational study comparing our approach with established benchmark compensation schemes from the literature. The results of our study show that the use of crowdshippers can lead to substantial benefits and that our approach outperforms the compensation schemes from the literature in terms. Our model allows operators to offer individualized compensation that results in a very high rate of accepted offers while achieving higher cost savings than the other schemes, which is due to the flexibility of our model according to the results. This is an important observation from the perspective of both the operator and the occasional driver. It shows that more occasional drivers can be utilized, which may result in lower dependency on third party logistic companies or a dedicated private fleet to fulfill the tasks. Moreover, a high acceptance rate of offers shows that the offers are more persuasive to occasional drivers, which signals higher driver satisfaction and therefore may result in a higher engagement and availability of occasional drivers in the long run. The sensitivity analysis shows that the proposed model is robust with regard to variations in the availability of occasional drivers and adapts very well to changes in the magnitude of the penalty incurred for rejected offers as well as to changes in the utility functions of occasional drivers. 
We also showed that the \probabbr can be used as a powerful heuristic in settings in which drivers and tasks appear dynamically.

Future research can be devoted to further studying dynamic versions of the problem, where tasks and available drivers are unknown in advance and become available over time. Given its efficiency, our approach
could be embedded in an approximate dynamic programming framework to take compensation decisions at each time step while
acceptance probabilities could be dynamically updated according to the earlier decisions taken by occasional drivers. 
Another possible avenue is to enrich the problem by considering additional features such as offering bundles of tasks to occasional drivers.

\ifTR
\bibliographystyle{plainnat}
\fi

\bibliography{ref}

\appendix







\section{Proofs\label{appendix:proofs}}

\appendixproof{\cref{th:compensation-values}}
{
	We conclude that the second statement of the theorem holds since compensation values $C_{ij}$ are irrelevant if task $i\in I$ is not offered to occasional driver $j\in J$, cf.\ objective function~\eqref{eq:obj}.
	Assume that task $i\in I$ is offered to occasional driver $j\in J$. We observe that an optimal compensation value $C_{ij}^*$ is equal to 
	\[
	C_{ij}^* = \argmin_{C_{ij}\ge 0} P_{ij}(C_{ij}) (C_{ij} - c'_i) + c_i'
	\]		
	since this value minimizes the (relevant terms of) objective function~\eqref{eq:minlp} and there are no dependencies on other tasks and drivers. The theorem follows because the last term of this equation is constant and can therefore be neglected. 
}

\appendixproof{\cref{prop:compenation-values:linear}}
{
	Recall that the optimal compensation values are calculated according to Equation~\eqref{eq:compensations:optimal} given in \cref{corr:compensations:optimal}. 
	We observe that the minimum of Equation~\eqref{eq:compensations:optimal} when plugging in probability function~\eqref{eq:prob:linear:proof} is attained for $C_{ij}\in ]0,c_i']$, in which case its value is less than or equal to zero, while it is equal to zero if $C_{ij}=0$. Similarly, we have $C_{ij}\le \frac{1-\alpha_{ij}}{\beta_{ij}}$ since $P_{ij}(\frac{1-\alpha_{ij}}{\beta_{ij}})=1$ and further increasing the compensation would increase the value of $P_{ij}(C_{ij})(C_{ij}-c_i')$.
	Thus, an optimal compensation value for task $i\in I$ when offered to occasional driver $j\in J$ must be a minimizer of $(\alpha_{ij} + \beta_{ij} C_{ij})(C_{ij}-c_i') = \beta_{ij} C_{ij}^2 + (\alpha_{ij} - \beta_{ij} c_i') C_{ij} - \alpha_{ij} c_i'$. By taking the first derivative (and since the second derivative is non-negative), the minimizer $C_{ij}^*$ is obtained as
	\[
	2 \beta_{ij} C_{ij}^* + \alpha_{ij} - \beta_{ij}c_i' = 0 \Leftrightarrow
	C_{ij}^* = \frac{c_i'}{2} - \frac{\alpha_{ij}}{2 \beta_{ij}}.
	\]
}

\appendixproof{\cref{prop:compenation-values:logistic}}
{
	We first observe that optimal compensation values lie in the intervals $]0,c_i']$ since Equation~\eqref{eq:compensations:optimal} evaluates to a non-positive value in this interval. Thus, an optimal compensation value for task $i\in I$ when offered to occasional driver $j\in J$ must be a minimizer of $\frac{C_{ij} - c_i'}{1 + e^{ - (\gamma_{ij} + \delta_{ij} C_{ij})}}$ whose first derivative 
	$\frac{e^{-(\gamma_{ij}+\delta_{ij} C_{ij})}\left( 1+(C_{ij}-c_i')\delta_{ij}+e^{\gamma_{ij} + \delta_{ij} C_{ij}}\right)}{\left(1+e^{-(\gamma_{ij} + \delta_{ij} C_{ij})}\right)^2}$
	is equal to zero if and only if the (strictly monotonically increasing) function $1+(C_{ij}-c_i')\delta_{ij}+e^{\gamma_{ij} + \delta_{ij} C_{ij}}$ is equal to zero.
	We use the following basic manipulations
	\begin{align*}
		1+(C_{ij}^*-c_i')\delta_{ij}+e^{\gamma_{ij} + \delta_{ij} C_{ij}^*}	 = 0 & \quad \Leftrightarrow \quad 
		e^{\gamma_{ij} + \delta_{ij} C_{ij}^*}  = c_i' \delta_{ij} - 1 - \delta_{ij} C_{ij}^* \quad  \Leftrightarrow \\
		\Leftrightarrow \quad  (c_i' \delta_{ij} - 1 - \delta_{ij} C_{ij}^*) e^{-\gamma_{ij} - \delta_{ij} C_{ij}^*} = 1 & \quad \Leftrightarrow \quad (c_i' \delta_{ij} - 1 - \delta_{ij} C_{ij}^*) e^{c_i' \delta_{ij} - 1 - \delta_{ij} C_{ij}^*} = e^{\gamma_{ij} + c_i'\delta_{ij} - 1}
	\end{align*}
	that hold for optimal compensation values $C_{ij}^*$ and observe that the last equation has the form $z e^z = b$ for $z=c_i' \delta_{ij} - 1 - \delta_{ij} C_{ij}^*$ and $b=e^{\gamma_{ij} + c_i'\delta_{ij} - 1}$. It can therefore be solved using the real-valued Lambert $W$ function $W(z e^z)=z$ and it follows that 
	\begin{align*}
		c_i' \delta_{ij} - 1 - \delta_{ij} C_{ij}^* & = W(e^{\gamma_{ij} + c_i'\delta_{ij} - 1}) & \quad \Leftrightarrow \quad 
		C_{ij}^* = -\frac{W(e^{\gamma_{ij} + c_i'\delta_{ij} - 1}) - c_i' \delta_{ij} + 1}{\delta_{ij}}.
	\end{align*}
}

\appendixproof{\cref{prop:np}}
{
	Obviously, the problem is in NP, since any assignment and compliance with the constraints can be verified in polynomial time (given there is only a polynomial number of constraints of type \eqref{eq:non-separable}). We show the proposition by reduction from the multidimensional knapsack problem. Let $A$ be an instance of the multidimensional knapsack problem with $d$ dimensions and items $I^K$ with value $v^K_{i}$ for $i^K \in I^K$. Let $C^K_\ell$ be the capacity of the knapsack in dimension $\ell = 1,\ldots,d$ and $w^K_{i\ell}$ be the weight of $i^K \in I^K$ in dimension $\ell = 1,\ldots,d$. Design an instance $A'$ of \probnonsepabbr with tasks $I$ as follows: 
	\begin{itemize}
		\item Let $\omega = \Pi_{i \in I^K} \min\left\{v_{i^K},2\right\}$ be a large constant.
		\item With each item $i^K \in I^K$ identify a task $i \in I$ with $c_i = c'_i = \omega v^K_i$. 
		\item Let $J$ be the set of occasional drivers with $|J| \ge |I|$. For each $j \in J$ set the parameters of the probability function such that $j$ accepts a task $i\in I$ for any $C_{ij} > 0$.
		\item For $\ell = 1,\ldots,d$ introduce the non-separability constraint
		\begin{align}
			\sum_{i\in I} \sum_{j\in J} (w^K_{i\ell} x_{ij}) \le C^K_\ell.	\label{eq:non-separable_np}
		\end{align}
	\end{itemize}
	Each constraint \eqref{eq:non-separable_np} in $A'$ models the knapsack constraint for dimension $\ell$ in $A$. Then, there exists a solution to $A$ with a value of at least $V$ if and only if there exists an assignment of tasks in $A'$ with an expected cost of at most $\tilde{C} = \left(\sum_{i^k \in I^K} v_{i^k} - V\right) \omega + 1$.
	\begin{itemize}
		\item[$\Rightarrow$] Let $i_1^K,\ldots,i_n^K$ be a feasible solution for $A$ with $V = \sum_1^n v^K_i$, which are identified by tasks $1,\ldots,n$ in $A'$. 	Take $n$ occasional drivers in $J$, w.l.o.g. $1,\ldots,n$ and offer a compensation of $C_{ij} = \frac{1}{n}$. The construction of the probability function ensures that the occasional drivers will accept these tasks (with a probability of one). Due to the construction of the non-separability constraints, and since the solution is feasible for the $A$, assigning the tasks to these occasional drivers is also feasible in $A'$. Since no other tasks are offered to the occasional drivers, the total expected cost is $\tilde{C}$. 
		\item[$\Leftarrow$] Consider a feasible solution to $A'$ with cost of $C' \leq \tilde{C}$ in which, w.l.o.g., tasks $1,\ldots,n$ are allocated to occasional drivers. Since the occasional drivers require a compensation strictly greater than zero to perform the task, it holds that 
		\[
		\sum_{j \in J}\sum_{i \in J} v^K_i x_{ij} \geq V
		\]
		since otherwise the total expected cost would be higher than $\tilde{C}$. Then, the corresponding items $i^K_1,\ldots,i^K_n$ also sum up to a value of at least $V$ in $A$. Since all non-separability constraints are adhered to in $A'$, this solution is also feasible for $A$.
	\end{itemize}
}

\appendixproof{\cref{prop:piecewise}}
{	
	We first recall that $x_{ij}=0$ implies $C_{ij}=0$ due to inequalities~\eqref{eq:minlp:forcing} and that $P_{ij}(0)=0$ holds by assumption. Thus, the nonlinear terms $P_{ij}(C_{ij}) x_{ij}$ of the objective function~\eqref{eq:minlp:obj} are equal to $P_{ij}(C_{ij})$. 	
	By first multiplying out \eqref{eq:minlp:obj} and then replacing  $P_{ij}(C_{ij}) x_{ij}$ by $P_{ij}(C_{ij})$ we obtain
	%
	\begin{align*}
		\sum_{i\in I} c_i y_i + \sum_{i\in I} \sum_{j\in J}  \left( P_{ij}(C_{ij}) C_{ij} + (1-P_{ij}(C_{ij})) c_i' \right) x_{ij} 
		= &\sum_{i\in I} c_i y_i + \sum_{i\in I} \sum_{j\in J}  \left( C_{ij} P_{ij}(C_{ij}) x_{ij} - c_i' P_{ij}(C_{ij}) x_{ij} + c_i' x_{ij} \right) = \\
		= & \sum_{i\in I} c_i y_i + \sum_{i\in I} \sum_{j\in J}  \left( P_{ij}(C_{ij})(C_{ij} - c_i') + c_i' x_{ij} \right) 	
	\end{align*}

	\cref{prop:piecewise} follows from substituting $f(C_{ij})= \sum_{i\in I} \sum_{j\in J} P_{ij}(C_{ij})(C_{ij} - c_i')$ and $g_{ij}(x_{ij})= c_i' x_{ij}$.

}

\appendixproof{\cref{prop:piecewise:linear}}
{
	Using $U_{ij}= \min\{c_i, \frac{1-\alpha_{ij}}{\beta_{ij}}\}$ discussed in \cref{sec:mnlp:linear}, the linear acceptance probability function~\eqref{eq:prob:linear} can be simplified to $P_{ij}(C_{ij})=\alpha_{ij} + \beta_{ij}C_{ij}$ due to constraints~\eqref{eq:minlp:forcing}, cf.\ \cref{prop:piecewise}. Thus, $f_{ij}(C_{ij}) + g_{ij}(x_{ij})$ must be equal to the corresponding (nonlinear) part of \eqref{eq:mnlp:obj:linear-prob} for each $i\in I$ and $j\in J$ and we have
	\begin{align*}
		f_{ij}(C_{ij}) + g_{ij}(x_{ij}) & = (P_{ij}(C_{ij}) C_{ij} + (1- P_{ij}(C_{ij})) c_i') x_{ij} = \\
		& = ((\alpha_{ij}+\beta_{ij} C_{ij}) C_{ij}  + (1-\alpha_{ij}-\beta_{ij} C_{ij}) c'_i) x_{ij} = \\
		& = \beta_{ij} C_{ij}^2 x_{ij} + (\alpha_{ij} - \beta_{ij} c_i' C_{ij} x_{ij} + c_i'(1 - \alpha_{ij}) x_{ij} = \\	
		& = \underbrace{\beta_{ij} C_{ij}^2 + (\alpha_{ij} - \beta_{ij} c_i') C_{ij}}_{f_{ij}(C_{ij})} + \underbrace{c_i'(1 - \alpha_{ij}) x_{ij}}_{g_{ij}(x_{ij})}
	\end{align*}
	
	The last equation holds since $x_{ij}=0$ implies that $C_{ij}=0$ due to constraints~\eqref{eq:minlp:forcing}. The second part of the proposition follows since $f_{ij}(C_{ij})$ is a convex (quadratic) function as $\beta_{ij}>0$.
}

\appendixproof{\cref{prop:piecewise:logistic}}
{
	We first observe that the proposition holds for $C_{ij}=0$ in which case we obtain $f_{ij}(C_{ij})+g_{ij}(x_{ij}) = c_i' x_{ij}$. For $0 < C_{ij}\le U_{ij}$ we obtain
	\begin{align*}	
		f_{ij}(C_{ij}) + g_{ij}(x_{ij}) & = (P_{ij}(C_{ij}) C_{ij} + (1- P_{ij}(C_{ij})) c_i') x_{ij} 
		= \frac{(C_{ij} - c_i') x_{ij}}{1+e^{-\gamma_{ij} - \delta_{ij} C_{ij}}} + c_i' x_{ij}.
	\end{align*}
	We observe that we can remove $x_{ij}$ from $(C_{ij} - c_i') x_{ij}$ since the assumption $C_{ij}> 0$ implies that $x_{ij}=1$. 
}




\section{Additional and detailed computational results\label{appendix:results}}

\begin{longtable}{rrlllrrrrrrr}
	\kill
	\caption{Paired t-test results comparing the individualized compensation scheme with the benchmark schemes in terms of total expected cost obtained using linear acceptance probability.} \\
	\toprule
	\multirow{2}{*}{$\odno$} & \multirow{2}{*}{$\pnlty$} & \multirow{2}{*}{$\utlty$} & \multirow{2}{*}{Individualized} & \multicolumn{2}{c}{Detour}& \multicolumn{2}{c}{Distance}& \multicolumn{2}{c}{Flat} \\ \cmidrule(lr){5-6} \cmidrule(lr){7-8} \cmidrule(lr){9-10}
	& & & & $\%$ Diff. & p-Val & $\%$ Diff. & p- Val & $\%$ Diff. & p-Val \\ 
	\midrule
	\endfirsthead
	\caption{(continued.)} \\
	\toprule
	\multirow{2}{*}{$\odno$} & \multirow{2}{*}{$\pnlty$} & \multirow{2}{*}{$\utlty$} & \multirow{2}{*}{Individualized} & \multicolumn{2}{c}{Detour}& \multicolumn{2}{c}{Distance}& \multicolumn{2}{c}{Flat} \\ \cmidrule(lr){5-6} \cmidrule(lr){7-8} \cmidrule(lr){9-10}
	& & & & $\%$ Diff. & p-Val & $\%$ Diff. & p- Val & $\%$ Diff. & p-Val \\ 
	\midrule
	\endhead
	50 & 0.00  & 0.00 & 4617.53 & 0.15 & 0.00 & 0.08 & 0.00 & 0.08 & 0.00 \\ 
	&  & 0.10 & 4573.57 & 0.14 & 0.00 & 0.08 & 0.00 & 0.08 & 0.00 \\ 
	&  & 0.20 & 4526.51 & 0.13 & 0.00 & 0.08 & 0.00 & 0.08 & 0.00 \\ 
	&  & 0.30 & 4476.96 & 0.12 & 0.00 & 0.08 & 0.00 & 0.08 & 0.00 \\ 
	&  & 0.40 & 4424.89 & 0.12 & 0.00 & 0.07 & 0.00 & 0.07 & 0.00 \\ 
	&  & 0.50 & 4370.55 & 0.11 & 0.00 & 0.07 & 0.00 & 0.07 & 0.00 \\ 
	&  & 0.60 & 4311.78 & 0.10 & 0.00 & 0.07 & 0.00 & 0.06 & 0.00 \\ 
	&  & 0.70 & 4248.36 & 0.09 & 0.00 & 0.07 & 0.00 & 0.06 & 0.00 \\
\end{longtable}

\begin{longtable}{rrlllrrrrrrr}
	\kill
	\caption{Paired t-test results comparing the individualized compensation scheme with the benchmark schemes in terms of mean acceptance rate obtained using linear acceptance probability.} \\
	\toprule
	\multirow{2}{*}{$\odno$} & \multirow{2}{*}{$\pnlty$} & \multirow{2}{*}{$\utlty$} & \multirow{2}{*}{Individualized} & \multicolumn{2}{c}{Detour}& \multicolumn{2}{c}{Distance}& \multicolumn{2}{c}{Flat} \\ \cmidrule(lr){5-6} \cmidrule(lr){7-8} \cmidrule(lr){9-10}
	& & & & $\%$ Diff. & p-Val & $\%$ Diff. & p- Val & $\%$ Diff. & p-Val \\ 
	\midrule
	\endfirsthead
	\caption{(continued.)} \\
	\toprule
	\multirow{2}{*}{$\odno$} & \multirow{2}{*}{$\pnlty$} & \multirow{2}{*}{$\utlty$} & \multirow{2}{*}{Individualized} & \multicolumn{2}{c}{Detour}& \multicolumn{2}{c}{Distance}& \multicolumn{2}{c}{Flat} \\ \cmidrule(lr){5-6} \cmidrule(lr){7-8} \cmidrule(lr){9-10}
	& & & & $\%$ Diff. & p-Val & $\%$ Diff. & p- Val & $\%$ Diff. & p-Val \\ 
	\midrule
	\endhead
	50 & 0.00 & 0.00 & 0.91 & -0.01 & 0.62 & -0.13 & 0.00 & -0.06 & 0.02 \\ 
	&      & 0.10 & 0.92 & -0.04 & 0.11 & -0.13 & 0.00 & -0.06 & 0.00 \\ 
	&      & 0.20 & 0.92 & -0.03 & 0.23 & -0.13 & 0.00 & -0.05 & 0.00 \\ 
	&      & 0.30 & 0.92 & -0.02 & 0.45 & -0.13 & 0.00 & -0.06 & 0.00 \\ 
	&      & 0.40 & 0.92 & -0.02 & 0.47 & -0.12 & 0.00 & -0.05 & 0.00 \\ 
	&      & 0.50 & 0.92 & -0.02 & 0.50 & -0.11 & 0.00 & -0.05 & 0.00 \\ 
	&      & 0.60 & 0.92 & -0.02 & 0.40 & -0.10 & 0.00 & -0.05 & 0.00 \\ 
	&      & 0.70 & 0.92 & -0.03 & 0.17 & -0.11 & 0.00 & -0.05 & 0.00 \\
\end{longtable}

\begin{longtable}{rrlllrrrrrrr}
	\kill
	\caption{Paired t-test results comparing the individualized compensation scheme with the benchmark schemes in terms of total expected cost obtained using logistic acceptance probability.} \\
	\toprule
	\multirow{2}{*}{$\odno$} & \multirow{2}{*}{$\pnlty$} & \multirow{2}{*}{$\utlty$} & \multirow{2}{*}{Individualized} & \multicolumn{2}{c}{Detour}& \multicolumn{2}{c}{Distance}& \multicolumn{2}{c}{Flat} \\ \cmidrule(lr){5-6} \cmidrule(lr){7-8} \cmidrule(lr){9-10}
	& & & & $\%$ Diff. & p-Val & $\%$ Diff. & p- Val & $\%$ Diff. & p-Val \\ 
	\midrule
	\endfirsthead
	\caption{(continued.)} \\
	\toprule
	\multirow{2}{*}{$\odno$} & \multirow{2}{*}{$\pnlty$} & \multirow{2}{*}{$\utlty$} & \multirow{2}{*}{Individualized} & \multicolumn{2}{c}{Detour}& \multicolumn{2}{c}{Distance}& \multicolumn{2}{c}{Flat} \\ \cmidrule(lr){5-6} \cmidrule(lr){7-8} \cmidrule(lr){9-10}
	& & & & $\%$ Diff. & p-Val & $\%$ Diff. & p- Val & $\%$ Diff. & p-Val \\ 
	\midrule
	\endhead
	50 & 0.00  & 0.00 & 6270.38 & 0.02 & 0.00 & 0.02 & 0.00 & 0.03 & 0.00 \\ 
	&  & 0.10 & 6182.13 & 0.02 & 0.00 & 0.02 & 0.00 & 0.03 & 0.00 \\ 
	&  & 0.20 & 6119.91 & 0.02 & 0.00 & 0.02 & 0.00 & 0.03 & 0.00 \\ 
	&  & 0.30 & 5980.57 & 0.01 & 0.00 & 0.02 & 0.00 & 0.03 & 0.00 \\ 
	&  & 0.40 & 5874.89 & 0.01 & 0.00 & 0.02 & 0.00 & 0.03 & 0.00 \\ 
	&  & 0.50 & 5775.77 & 0.01 & 0.00 & 0.02 & 0.00 & 0.03 & 0.00 \\ 
	&  & 0.60 & 5666.36 & 0.01 & 0.00 & 0.02 & 0.00 & 0.03 & 0.00 \\ 
	&  & 0.70 & 5545.61 & 0.01 & 0.00 & 0.02 & 0.00 & 0.03 & 0.00 \\
\end{longtable}

\begin{longtable}{rrlllrrrrrrr}
	\kill
	\caption{Paired t-test results comparing the individualized compensation scheme with the benchmark schemes in terms of mean acceptance rate obtained using logistic acceptance probability.} \\
	\toprule
	\multirow{2}{*}{$\odno$} & \multirow{2}{*}{$\pnlty$} & \multirow{2}{*}{$\utlty$} & \multirow{2}{*}{Individualized} & \multicolumn{2}{c}{Detour}& \multicolumn{2}{c}{Distance}& \multicolumn{2}{c}{Flat} \\ \cmidrule(lr){5-6} \cmidrule(lr){7-8} \cmidrule(lr){9-10}
	& & & & $\%$ Diff. & p-Val & $\%$ Diff. & p- Val & $\%$ Diff. & p-Val \\ 
	\midrule
	\endfirsthead
	\caption{(continued.)} \\
	\toprule
	\multirow{2}{*}{$\odno$} & \multirow{2}{*}{$\pnlty$} & \multirow{2}{*}{$\utlty$} & \multirow{2}{*}{Individualized} & \multicolumn{2}{c}{Detour}& \multicolumn{2}{c}{Distance}& \multicolumn{2}{c}{Flat} \\ \cmidrule(lr){5-6} \cmidrule(lr){7-8} \cmidrule(lr){9-10}
	& & & & $\%$ Diff. & p-Val & $\%$ Diff. & p- Val & $\%$ Diff. & p-Val \\ 
	\midrule
	\endhead
	50 & 0.00 & 0.00 & 0.51 & -0.07 & 0.01 & -0.12 & 0.00 & -0.18 & 0.00 \\ 
	&      & 0.10 & 0.53 & -0.06 & 0.02 & -0.11 & 0.00 & -0.18 & 0.00 \\ 
	&      & 0.20 & 0.54 & -0.05 & 0.02 & -0.11 & 0.00 & -0.17 & 0.00 \\ 
	&      & 0.30 & 0.56 & -0.05 & 0.02 & -0.10 & 0.00 & -0.15 & 0.00 \\ 
	&      & 0.40 & 0.58 & -0.04 & 0.02 & -0.10 & 0.00 & -0.14 & 0.00 \\ 
	&      & 0.50 & 0.59 & -0.03 & 0.02 & -0.09 & 0.00 & -0.12 & 0.00 \\ 
	&      & 0.60 & 0.61 & -0.01 & 0.25 & -0.09 & 0.00 & -0.12 & 0.00 \\ 
	&      & 0.70 & 0.63 & -0.01 & 0.19 & -0.10 & 0.00 & -0.10 & 0.00 \\
\end{longtable}

\end{document}